\documentclass[onefignum,onetabnum]{siamonline190516}

\usepackage{lipsum}
\usepackage{amsfonts}
\usepackage{ amssymb }
\usepackage{graphicx}
\usepackage{epstopdf}
\usepackage{algorithmic}
\usepackage{subcaption}

\ifpdf
  \DeclareGraphicsExtensions{.eps,.pdf,.png,.jpg}
\else
  \DeclareGraphicsExtensions{.eps}
\fi

\usepackage{enumitem}
\setlist[enumerate]{leftmargin=.5in}
\setlist[itemize]{leftmargin=.5in}

\newcommand{\rd}{\mathrm{d}}

\newcommand{\ql}[1]{{\color{orange}{#1}}}

\newsiamremark{remark}{Remark}
\newsiamremark{hypothesis}{Hypothesis}
\crefname{hypothesis}{Hypothesis}{Hypotheses}
\newsiamthm{claim}{Claim}

\headers{Reconstructing the phonon transmission coefficient at solid interfaces}{}

\title{Reconstructing the thermal phonon transmission coefficient at solid interfaces in the phonon transport equation\thanks{Submitted to the editors DATE.
\funding{The research of I.M.G is supported by NSF DMS 2009736 and RNMS 1107465. The research of Q.L and A.N is partially supported by NSF DMS 1750488, ONR-N00014-21-1-2140 and RNMS 1107465. All authors thank and gratefully acknowledge the hospitality and support from the Oden Institute of Computational Engineering and Sciences at the University of Texas Austin. All authors thank the two anonymous reviewers for insightful suggestions.}}}
\author{Irene M. Gamba\thanks{Department of Mathematics and Oden Institute, University of Texas-Austin, Austin, TX 78712 (gamba@math.utexas.edu) }
\and Qin Li\thanks{Department of Mathematics, University of Wisconsin-Madison, Madison, WI 53706 (qinli@math.wisc.edu)}
\and Anjali Nair\thanks{Department of Mathematics, University of Wisconsin-Madison, Madison, WI 53706 (nair25@wisc.edu) }}

\usepackage{amsopn}

\ifpdf
\hypersetup{
  pdftitle={Reconstructing thermal phonon transmission coefficient at solid interfaces},
  pdfauthor={Irene Gamba, Qin Li, Anjali Nair}
}
\fi

\begin{document}

\maketitle

\begin{abstract}

The ab initio model for heat propagation is the phonon transport equation, a Boltzmann-like kinetic equation. When two materials are put side by side, the heat that propagates from one material to the other experiences thermal boundary resistance. Mathematically, it is represented by the reflection coefficient of the phonon transport equation on the interface of the two materials. This coefficient takes different values at different phonon frequencies, between different materials. In experiments scientists measure the surface temperature of one material to infer the reflection coefficient as a function of phonon frequency.

In this article, we formulate this inverse problem in an optimization framework and apply the stochastic gradient descent (SGD) method for finding the optimal solution. We furthermore prove the maximum principle and show the Lipschitz continuity of the Fr\'echet derivative. These properties allow us to justify the application of SGD in this setup.
\end{abstract}



\section{Introduction}

How heat propagates is a classical topic. Mathematically a parabolic type heat equation has long been regarded as the model to describe heat conductance. It was recently discovered that the underlying ab initio model should be the phonon transport Boltzmann equation~\cite{hua2017experimental}. One can formally derive that this phonon transport equation degenerates to the heat equation in the macro-scale regime when one assumes that the Fourier law holds true. This is to assume that the rate of heat conductance through a material is negatively proportional to the gradient of the temperature.

Rigorously deriving the heat equation from kinetic equations such as the linear Boltzmann equation or the radiative transfer equation is a standard process in the so-called parabolic regime~\cite{bardos1984diffusion}. That means for light propagation (using the radiative transfer equation), we can link the diffusion effect with the light scattering. However, such derivation does not exist for the phonon system. One problem we encounter comes from the fact that phonons, unlike photons, have all frequencies $\omega$ coupled up. In particular, the collision operator for the phonon transport equation contributes an equilibrium term that satisfies a Bose-Einstein distribution in the frequency domain: it allows energy contained in one frequency to transfer to another. Furthermore the speed in the phonon transport term involves not only the velocity $\mu$, but also the group velocity $v(\omega)$ which also positively depends on the frequency $\omega$. These differences make the derivation of the diffusion limit (or the parabolic derivation) for phonon transport equation not as straightforward.

On the formal level, passing from the ab initio Boltzmann model to the heat equation requires the validity of the Fourier law, an assumption that breaks down when kinetic phenomenon dominates. This is especially true at interfaces where different solid materials meet. Material discontinuities lead to thermal phonon reflections. On the macroscopic level, it is observed as thermal boundary resistance (TBR), and is reflected by a temperature drop at the interface~\cite{cahill2003nanoscale}. TBR exists at the interface between two dissimilar materials due to differences in phonon states on each side of the interface. Defects and roughness further enhance phonon reflections and TBR effect. Such effect can hardly be explained, or measured directly on the heat equation level.

As scientists came to the realization that the underlying ab initio model is the Boltzmann equation instead of the heat equation, more and more experimental studies are conducted to reveal the model parameters and properties of the Boltzmann equation. In the recent years, a lot of experimental work has been done to understand the heat conductance inside a material or at the interface of two solids, hoping these collected data can help in designing materials that have better heat conductance or certain desired heat properties~\cite{lyeo2006thermal, norris2009examining, wang2011thermal}.


\begin{figure}[htb]
\centering
\includegraphics[scale = 0.75]{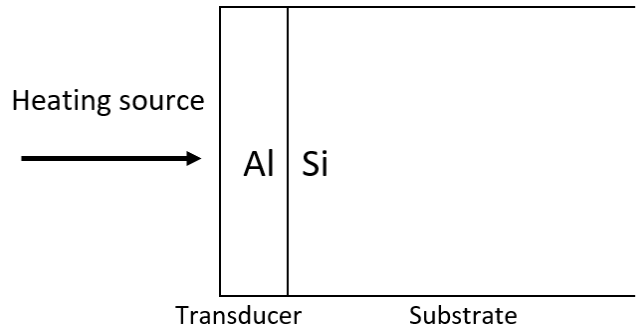}
\caption{Experiment setup: details can be found in ~\cite{hua2017experimental}. In experiments, two solid materials are placed side by side, and heat is injected on the surface of aluminium. Temperature is also measured at the same location as a function of time.}
\label{fig:illustration}
\end{figure}

This is a classical inverse problem. Measurements are taken to infer the material properties. In~\cite{hua2017experimental} the authors essentially used techniques similar to least square with an $L_1$ penalty term to ensure sparsity. As the first investigation into this problem, the approach gives a rough estimate of the parameter, but mathematically it is very difficult to evaluate the accuracy of the recovery. In this article, we study this problem with a more rigorous mathematical view. We will also confine ourselves in the optimization framework. Given certain amount of data (certain number of experiments, and certain number of measurements in one experiment), we formulate the problem into a least square setup to have the mismatch minimized. We will study if this minimization problem is well-posed, and how to numerically find the solution in an efficient way. More specifically, we will first formulate the optimization problem, apply the stochastic gradient descent method, the state-of-the-art optimization technique, and then derive the problem's Fr\'echet derivative. This allows us to investigate the associated convergence property. We will demonstrate that the system has maximum principle, which allows us to justify the objective function to be convex in a small region around the optimal point. This further means the SGD approach will converge. Finally we apply the method we derived in two examples in the numerical section.

We should note that we took a practical route in this paper and focus on the numerical property of the problem. Investigating the well-posedness of the inverse problem, such as proving the uniqueness and the stability of the reconstruction require much more delicate PDE analysis, and is left for future investigation.

{We should also note that there are quite some results on inverse problems for kinetic equations. The research is mostly focused on the reconstruction of optical parameters in the radiative transfer equation (RTE). In~\cite{CS96}, the authors pioneered the problem assuming the entire albedo operator is known. The stability was later proved in ~\cite{Romanov96,Wang,bal2010stability,LaiLiUhlmann} in different scenarios. In~\cite{bal2008inverse,time_harmonic,BalMonard_time_harmonic} the authors revised the measurement-taking mechanism and assumed only the intensity is known. They furthermore studied the time-harmonic perturbation to stablize the inverse problem. The connection of optical tomography and the Calder\'on problem was drawn in~\cite{chen2018stability, LaiLiUhlmann} where the stability and its dependence on the Knudsen number is investigated. See reviews in~\cite{Ren_review,arridge1999optical}. There are many other imaging techniques that involve the coupling of the radiative transfer equation with other equations, that lead to bioluminescence tomography, quantitative photoacoustic tomography, and acousto-optic imaging~\cite{chung2020coherent,ChungSchotland,Ren_fPAT,BalChungSchotland,BalSchotland2010,BalTamasan}, among many others. Outside the realm of the radiative transfer equation, the inverse problem is also explored for the Vlasov-Poisson-Boltzmann system that characterizes the dynamics of electrons in semi-conductors~\cite{Gamba}. The authors of the paper investigated how to use Neumann data and the albedo operator at the surface of semi-conductors to infer the deterioration in the interior of the semi-conductor. Recently, in~\cite{CarlemanEstimate} the authors used the Carleman estimate to infer the parameters in general transport type equation in the kinetic framework. We stress that these works investigate inverse problems for the kinetic equations, and the to-be-reconstructed parameters are usually related to the heterogeneous media. Photon frequency domain is typically untouched. This is different from the setup in our paper, where the media is homogeneous in $x$ since the material on both sides of the interface are pure enough in the lab experiments. However, the reflection index is a function of phonon frequency, and this dimension of heterogeneity should not be discounted as it was used to in the literature for RTE.}

There is also abundant research in numerical inverse problems. For kinetic equations, multiple numerical strategies have been investigated, including both gradient based method and hessian based method~\cite{arridge1999optical, Arridge_couple,dehghani1999photon}. For general PDE-based inverse problems, SGD has been a very popular technique for finding the minimizer. This is because most PDE-based inverse problems eventually are formulated as minimization problem with the loss function being in the format of a summation of multiple smaller loss functions. This is the exactly the format where SGD outperforms other optimization methods~\cite{chen2018online}.

The paper is organized as follows. In Section~\ref{sec:set_up}, we present the model equation and its linearization (around the room temperature). The formulation of the inverse problem is also described in this section. Numerically to solve this optimization problem, we choose to use the stochastic gradient descent method, which would require the computation of the Fr\'echet derivative in each iteration. In Section~\ref{sec:SGD} we discuss the self-adjoint property of the collision operator and derive the Fr\'echet derivative. This allows us to summarize the application of SGD at the end of this section. We discuss properties of the optimization formulation and the use of SGD in our setting in Section~\ref{sec:properties}. In particular, we will discuss the maximum principle in Section~\ref{sec:property_MP}, and Section~\ref{sec:property_SGD} is dedicated to the properties of SGD applied in this particular problem. In particular we will show the Lipschitz continuity of the Fr\'echet derivative. We conduct two sets of numerical experiments: first, assuming that the reflection coefficient is parametrized by a finite set of variables and second, assuming the reflection coefficient is a simple smooth function without any prior parametrization. SGD gives good convergence in both cases, and the results are shown in Section~\ref{sec:numerics}.

\section{Model and the inverse problem setup}\label{sec:set_up}
In this section we present the phonon transport equation and its linearization. We also set up the inverse problem as a PDE-constrained minimization problem.

\subsection{Phonon transport equation and linearization}

The phonon transport equation can be categorized as a BGK-type kinetic equation. { Denote by $F(t,x,\mu,\omega)$ the deviational distribution function of phonons (viewed as particles in this paper) at time $t\in\mathbb{R}^+$ on the phase space at location $x\in\mathbb{R}^n$, transport speed $\mu\in\mathbb{S}^{n-1}$ and frequency $\omega\in\mathbb{R}^+$.} Since $\mu$ has a constant amplitude, we normalize it and set it to be $1$. In labs the experiments are set to be plane-symmetric, and the problem becomes pseudo 1D, meaning $n=1$, and $\mu = \cos\theta\in[-1,1]$. The equation then writes as:

\begin{equation}\label{eqn:bgk_transport}
\partial_tF+\mu v\partial_xF=\frac{F^\ast-F}{\tau(\omega)}\,.
\end{equation}

The two terms on the left describe the phonon transport with velocity $\mu$ and group velocity $v=v(\omega)$. The term on the right is called the collision term, with $\tau=\tau(\omega)$ representing the relaxation time. It can be viewed equivalently to the BGK (Bhatnagar-Gross-Krook) operator in the classical Boltzmann equation for the rarefied gas. $F^\ast$ is the Bose-Einstein distribution, also termed the Maxwellian, and is defined by:
\begin{equation*}
F^\ast(t,x,\mu,\omega)=\frac{\hbar\omega D(\omega)}{e^{\frac{\hbar\omega}{ kT(t,x)}}-1}\,.
\end{equation*}
In the formula, $k$ is the Boltzmann constant, $\hbar$ is the reduced Planck constant, and $D(\omega)$ is the phonon density of states. The profile is a constant in $\mu$ and approximately exponential in $\omega$ (for $\omega$ big enough), and the rate of the exponential decay is uniquely determined by temperature $T(t,x)$. {$F$ is related to the phonon distribution function $f$ as $F=\hbar\omega D(\omega)f$ ~\cite{hao2009frequency}.}

{The system preserves energy, meaning the right hand side vanishes when the zeroth moment is being taken.} This uniquely determines the temperature of the system, namely:

\[
\int^\infty_{0}\int_{-1}^1\frac{F}{\tau}\rd{\mu}\rd\omega = \int^\infty_{0} \int_{-1}^1 \frac{F^\ast(t,x,\mu,\omega)}{\tau}\rd\mu\rd\omega\,.
\]

{In experiments, the temperature is typically kept around the room temperature and the variation is rather small~\cite{hua2017experimental}. This allows us to linearize the system around the room temperature.} Denote the room temperature to be $T_0$, and the associated Maxwellian is denoted 
\[
F^\ast_0 = \frac{\hbar\omega D(\omega)}{e^{\frac{\hbar\omega}{kT_0}}-1}\,.
\]

We linearize equation~\eqref{eqn:bgk_transport} by subtracting $\frac{F^\ast_0}{\tau}$ and adding it back, and call  $g = F-F^\ast_0$, then with straightforward calculation, since $F_0^\ast$ has no dependence on $t$ and $x$:

\begin{equation}\label{eqn:linear_PTE1}
\partial_tg+\mu v(\omega)\partial_xg=\frac{-g}{\tau} +\frac{1}{\tau}[F^\ast-F^\ast_0]\,.
\end{equation}
For $F^\ast$ and $F^\ast_0$ evaluated at very similar $T$ and $T_0$, we approximate
\[
F^\ast-F^\ast_0\approx \partial_{T}F^\ast(T_0)\Delta T\,,\quad\text{with}\quad \Delta T = T-T_0\,,
\]
where higher orders at $\mathcal{O}(\Delta T^2)$ are eliminated. Noting that
\begin{align}
\partial_TF^\ast(T_0) =\frac{(\hbar\omega)^2D(\omega)}{kT_0^2}\frac{e^{\hbar\omega/kT_0}}{\left(e^{\hbar\omega/kT_0}-1\right)^2}\triangleq g^\ast
\end{align}
we rewrite equation~\eqref{eqn:linear_PTE1} to:
\begin{equation}\label{eqn:linear_PTE21}
\partial_tg+\mu v(\omega)\partial_xg=\frac{-g}{\tau}+\frac{g^\ast}{\tau}\Delta T\,.
\end{equation}
Similar to the nonlinear case, conservation law requires the zeroth moment of the right hand side to be zero, which amounts to 
\[
\int \frac{g^\ast}{\tau}(\omega)\rd\mu\rd\omega\Delta T =\int \frac{g}{\tau}(\omega)\rd\mu\rd\omega\quad\Rightarrow\quad \Delta T = \frac{\langle g/\tau\rangle}{\langle g^\ast/\tau\rangle}\,,
\]
where we use the bracket notation $\langle\rangle=\int_{\mu=-1}^1\int^\infty_{\omega=0}\rd\omega \rd\mu$.

For simplicity, we non-dimensionalize the system, and set $\hbar/kT_0 = 1$, and $D= 1$. We also make the approximation that $v(\omega) = \omega$ and $\tau\approx \frac{1}{\omega}$, as suggested in~\cite{PRB} and ~\cite{hua2017experimental}. We now finally arrive at
\begin{equation}\label{eqn:linear_PTE2}
 \partial_tg+\mu \omega\partial_xg=-\omega g + M\langle \omega g\rangle\,,
\end{equation}
with
\begin{equation}\label{eqn:L}
M = \frac{\omega g^\ast}{\langle\omega g^\ast\rangle}\,,\quad\text{and}\quad\Delta T = \frac{\langle \omega g\rangle}{\langle \omega g^\ast\rangle}\quad\text{and}\quad g^\ast=\frac{\omega^2 e^\omega}{(e^\omega-1)^2}\,.
\end{equation}
Clearly, the Maxwellian $M$ is normalized: $\langle M\rangle=1$. We also call the right hand side the linearized BGK operator
\begin{equation}\label{eqn:collision}
\mathcal{L} g=\frac{-g}{\tau}+\frac{\langle g/\tau\rangle}{\langle g^\ast/\tau\rangle}g^\ast(\omega)=-\omega g + M\langle \omega g\rangle\,.
\end{equation}

\begin{figure}
\centering
\includegraphics[scale = 0.4]{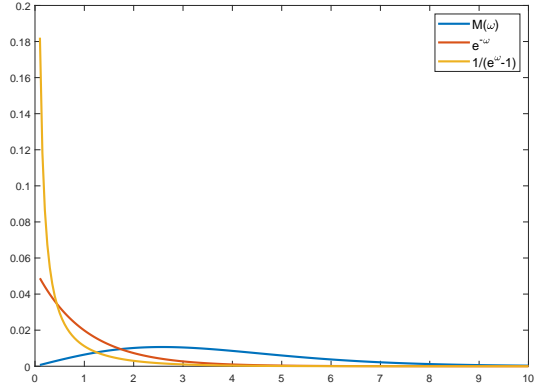} 
\caption{
{Plot of the Maxwellian $M$ as a function of $\omega$ (normalized). We compare it to the typically used Bose-Einstein distribution $1/(e^\omega-1)$ (normalized)and its approximation $e^{-\omega}$ (normalized) for quantum systems.} }
\label{fig:equilibrium}
\end{figure}

We note that the Maxwellian here is not a traditional one. {In Figure~\ref{fig:equilibrium} we depict the Maxwellian, and its comparison to some other standard Maxwellian functions.}

In practice, to handle the situation in Figure~\ref{fig:illustration}, phonon is in charge of heat transfer and the model equation is used in both the aluminum and silicon regions, i.e, we write $g$ and $h$ to be the density of phonon in the two regions respectively. Considering the boundary condition for the aluminum component is incoming-type and the boundary condition at the interface is reflective type, we have the following model (denote the location of the surface of aluminum to be $x=0$ and the interface to be $x=1$):
\begin{equation}\label{eqn:linear_PTE}
\begin{aligned}
\begin{cases}
\partial_tg+\mu v(\omega)\partial_xg=\mathcal{L}g\,,\quad x\in[0,1]\\
g(x=0\,,\mu\,,\cdot) = {\phi}\,,\quad \mu>0\\
g(x=1\,,\mu\,,\cdot)=\eta(\omega)g(x=1\,,-\mu\,,\cdot)\,,\quad \mu<0
\end{cases}\,,
\\
\begin{cases}
\partial_th+\mu v(\omega)\partial_xh=\mathcal{L}h\,,\quad x\in[1,\infty)\\
h(x=1\,,\mu\,,\cdot) = (1-\eta)g(x=1\,,\mu\,,\cdot)\,,\quad \mu>0
\end{cases}\,.
\end{aligned}
\end{equation}

In the model $g$ and $h$ are phonon density distribution within aluminum and silicon respectively. At $x=0$, laser beam injects heat into aluminum, and that serves as incoming boundary condition at $g(t, x=0,\mu>0,\omega)$, which we name by $\phi$. Due to the large size of silicon, it can be viewed that $h$ equation is supported on the half-domain. At $x=1$, the interface, we model the reflection coefficient to be $\eta$, meaning $\eta$ portion of the phonon density is reflected back with the flipped sign $\mu\to -\mu$. According to~\cite{hua2017experimental}, this reflection coefficient only depends on $\omega$, the frequency. It is the coefficient that determines how much heat gets propagated into silicon, and is the coefficient we need to reconstruct in the lab using measurements. {In lab experiments, the materials are large enough for the pseudo-1D assumption to hold true. This is seen in discussion in~\cite{hua2017experimental}. In reality, two kinds of materials can certainly touch each other through a curved interface and the above mentioned model no longer holds. However, the reflection index is only a function of frequency $\omega$ and thus the reconstructed $\eta(\omega)$ using this model is still valid.}

\subsection{Formulating the inverse problem}
The experiments are non-intrusive, in the sense that the materials are untouched, and the data is collected only on the surface $x=0$. The data we can tune is the incoming boundary condition $\phi$ in~\eqref{eqn:linear_PTE} and the data we can collect is the temperature at the surface, namely $\Delta T =\frac{\langle \omega g\rangle}{\langle \omega g^\ast\rangle}$ as a function of $t$ and $x=0$. In labs we can send many different profiles of $\phi$ into the system, and measure $\Delta T$ for the different $\phi$.

Accordingly, the forward map is:
\begin{equation}\label{map}
\mathcal{M}_\eta: {\phi} \to \Delta T(x=0,t)\,.
\end{equation}
It maps the incoming data configuration $\phi$ to the temperature at the surface as a function of time. The subscript $\eta$ represents the dependence of the map. In the forward setting, $\eta$ is assumed to be known, then equation~\eqref{eqn:linear_PTE} can be solved for any given $\phi$ for the outcome $\Delta T$. In reality, $\eta$ is unknown, and we test the system with a number of different configurations of ${\phi}$ and measure the corresponding $\Delta T$ to infer $\eta$. {We stress that the inverse problem is conducted on the linearized equation~\eqref{eqn:linear_PTE}. The laser that sends energy into the solid materials is typically not powerful enough to increase the temperature of the solid drastically away from the room temperature, and linearizing the system around the room temperature is a valid assumption~\cite{hua2017experimental}.}

While the to-be-reconstructed $\eta$ is a function of $\omega$, chosen in $L_\infty$, an infinite dimensional function space, there are infinitely many configurations of $\phi$ too. At this moment whether tuning $\phi$ gives a unique reconstruction of $\eta$ is an unknown well-posedness problem that we plan to investigate in the near future. In the current paper we mainly focus on the numerical and practical setting, namely, suppose a finite number of experiments are conducted with finitely many configurations of $\phi$ imposed, and in each experiment, the measurement $\Delta T$ is taken on discrete time, how do we reconstruct $\eta$?

We denote by $I$ the number of configurations of the boundary condition in different experiments, namely
\[
\{{\phi}_i\,,\quad i=1\,,\cdots,I\}
\]
is injected into~\eqref{eqn:linear_PTE} at different rounds of experiments. We also take measurements on the surface $x=0$ through test functions. Denote $\{\psi_j(t)\}$ the test function, the data points are then:
\[
\int \Delta T(t,x=0) \psi_j(t)\rd{t}\,,\quad\text{with}\quad j=1\,,\cdots,J\,.
\]
When we choose $\psi_j = \delta_{t_j}$ the temperature is simply taken at discrete time $\{t_j\}$.

We denote the data $d_{ij}$ to be $\Delta T$, collected at the $i$-th experiment, measured with test function $\psi_j$ with additive noise, then, for indices $(i,j)\in[0,I]\times[0,J]$
\[
d_{ij} = \mathcal{M}_{ij}(\eta) +\text{noise}\,,\quad\text{with}\quad \mathcal{M}_{ij}(\eta) = \int \Delta T_i(x=0,t)\psi_j(t)\rd{t}\,,
\]
where $\Delta T_i$ is the solution to~\eqref{eqn:linear_PTE} equipped with boundary condition ${\phi}_i$, assuming the reflection coefficient is $\eta$. The noise term is inevitable in experiments.

A standard approach to reconstruct $\eta(\omega)$ is to formulate a PDE-constrained minimization problem. We look for a function $\eta$ that minimizes the mismatch between the produced $\Delta T$ and the measured data:

\begin{equation}\label{eqn:opt_1}
\begin{aligned}
&\min_{\eta}\frac{1}{IJ}\sum_{ij}|\int \Delta T_i(x=0,t)\psi_j\rd{t}-d_{ij}|^2\\
\text{s.t.}\quad &
\begin{cases}
\partial_t g_i + \mu\omega\partial_x g_i = \omega g_i^\ast\Delta T_i-\omega g_i\,,\quad\text{with}\quad \Delta T_i (x,t)=\frac{\langle\omega g_i\rangle}{\langle \omega g^\ast\rangle}(x,t)\\
g_i(x=0,\mu,\cdot) = {\phi}_i\quad \mu>0\\
g_i(x=1,\mu,\cdot) = \eta g_i(x=1,-\mu,\cdot)\quad\mu<0
\end{cases}\,.
\end{aligned}
\end{equation}
Here $\frac{1}{IJ}$ is merely a constant and does not affect the minimum configuration, but adding it puts the formulation in the framework that SGD deals with.

In a concise form, define the loss function
\begin{equation}\label{eqn:def_L}
L = \frac{1}{IJ}\sum_{ij}L_{ij}^2\,,\quad\text{and}\quad L_{ij} = \mathcal{M}_{ij}(\eta)-d_{ij}\,,
\end{equation}
then the PDE-constrained minimization problem can also be written as:
\begin{equation}\label{eqn:optimization}
\min_{\eta}\frac{1}{IJ}\sum_{ij}|\mathcal{M}_{ij}(\eta)-d_{ij}|^2=\min_{\eta(\omega)}\frac{1}{IJ}\sum_{ij}L_{ij}^2=\min_{\eta(\omega)}L\,.
\end{equation}

We denote the minimizer to be $\eta_\ast$.

\begin{remark}
We have a few comments:
\begin{itemize}
\item If some prior information is known about $\eta$, this information could be built into the minimization formulation as a relaxation term. For example, if it is known ahead of time that $\eta$ should be close to $\eta_0$ in some sense, then the minimization is modified to
\begin{equation*}
\min_{\eta}\frac{1}{IJ}\sum_{ij}|\mathcal{M}_{ij}(\eta)-d_{ij}|^2+\lambda\|\eta-\eta_0\|^p_p\,.
\end{equation*}
where $L_p$ norm is chosen according to properties from physics and $\lambda$ is the relaxation coefficient. We assume in this paper that we do not have such prior information.
\item {In~\cite{hua2017experimental}, the experimentalists chose to model the system as a pseudo-1D system with plane geometry. The system would be modified if the geometry is changed to accommodate non-trivial curvature in high dimensions. This not only brings mathematical difficulty, but also makes the lab experiments much harder. One needs to deal with the artificial difficulties induced by the inhomogeneous spatial dependence.}
\item {We study the inverse problem where the forward model is linearized around the room temperature. This is a valid assumption since the energy injected into the system is typically not strong enough to trigger high temperature fluctuation. The inverse problem, however, evaluates $\mathcal{M}$'s dependence on $\eta$, and nevertheless is still nonlinear.}
\end{itemize}
\end{remark}

\section{Stochastic gradient descent}\label{sec:SGD}
There are many approaches for solving PDE-constrained minimization problems such as~\eqref{eqn:opt_1}, or equivalently~\eqref{eqn:optimization}. Most of the approaches are either gradient-based or Hessian-based. While usually gradient-based methods converge at a slower rate than methods that incorporate Hessians, the cost of evaluating Hessians is very high. For PDE-constrained minimization problems, every data point provides a different Hessian term, and there are many data points. Furthermore, if the to-be-reconstructed parameter is a function, the Hessian is infinite-dimensional, which makes a very large sized matrix upon discretization. This cost is beyond what a typical computer can afford. Thus, we choose gradient-based methods.

Among all gradient-based optimization methods, the stochastic gradient descent (SGD) started gaining ground in the past decade. It is a method that originated in the 90s~\cite{leon1998online} (which also sees its history back in~\cite{robbins1951stochastic}), and gradually became popular in the new data science era. As a typical example of probabilistic-type algorithms, it sacrifices certain level of accuracy in trade of efficiency. Unlike the standard Gradient Descent method, SGD does require a certain form of the objective function. $L(\eta)$ is an average of many smaller objective functions, meaning $L = \frac{1}{N}\sum_{i=1}^NL_i(\eta)$ with a fairly large $N$.

If a classical gradient descent method is used, then per iteration, to update $\eta_{n+1}$ from $\eta_n$, one needs to compute the gradient of $L$ at this specific $\eta_n$, and that amounts to computing $N$ gradients: $\nabla L_1\,,\cdots\nabla L_N$. For a large $N$, the cost is very high. For PDE-constrained minimization in particular, this $\nabla L_i$ is a Fr\'echet derivative, and is usually obtained through two PDE solves: one for the forward problem and one for the adjoint. $N$ gradients means $2N$ PDE-solves at each iteration. This cost is prohibitively high.

The idea of SGD is that in each iteration, only one $L_i$ is randomly chosen as a fair representative of $L$ and the single gradient $\nabla L_i$ is viewed as a surrogate of the entire $\nabla L$. So per iteration instead of computing $N$ gradients, one only computes $1$. This significantly reduces the cost for each iteration. However, replacing $\nabla L$ by a random representative $\nabla L_i$ brings extra error, and as a consequence, a larger number of iterations is required. This leads to a delicate error analysis. In the past few years, in what scenario does SGD win over the traditional GD has been a popular topic, and many variations and extensions of SGD have also been proposed~\cite{needell2014stochastic, le2012stochastic, zhao2015stochastic}. Despite the rich literature for SGD as a stand-alone algorithm, its performance in the PDE-constrained minimization setting is mostly in lack. We only refer the readers to the review papers in general settings~\cite{bottou2010large, bottou2018optimization}.

We apply SGD method to our problem. That is to say, per each time step, we randomly choose one set of data pair and use it to adjust the evaluation of $\eta$. More specifically, at time step $n$, instead of using all loss functions, one selects a multi-index $\gamma_n=(i_n\,,j_n)$ at random and performs gradient descent determined by this particular data pair:
\begin{equation}\label{eqn:SGD}
\eta_{n+1} = \eta_n -2\alpha_n L_{\gamma_n}(\eta_n)\nabla_{\eta} L_{\gamma_n}(\eta_n)\,.
\end{equation}
Here $\alpha_n$ is the time-step to be adjusted according to the user's preference.

In the formula of~\eqref{eqn:SGD}, the evaluation of $L_{\gamma_n}(\eta_n)$ is straightforward. According to equation~\eqref{eqn:def_L}, it amounts to setting $\eta$ in~\eqref{eqn:linear_PTE} as $\eta_n$ and computing it with incoming data ${\phi}_{i_n}$ as the boundary condition, and test the solution at $x=0$ with $\psi_{j_n}$. How to compute $\nabla_\eta L_{\gamma_n}$, however, is less clear. Considering $d_{\gamma_n}$ is a given data point and has no dependence on $\eta$, it is the Fr\'echet derivative of the map $\mathcal{M}$ on $\eta$ at this particular $\gamma_n$:
\begin{equation}
    \begin{aligned}
    \nabla_{\eta} L_{\gamma_n}&=\nabla_\eta \mathcal{M}_{\gamma_n}(\eta)=\lim_{|\delta\eta|\to0}\frac{\mathcal{M}_{\gamma_n}(\eta+\delta \eta)-\mathcal{M}_{\gamma_n}(\eta)}{\delta \eta}\\
    &=\frac{\delta}{\delta \eta} \int\Delta T_{i_n}(x=0,t)\psi_{j_n}(t)\rd{t}\,.
    \end{aligned}
\end{equation}

The computation of this derivative requires the computation of the adjoint equation. Before stating the result, we first notice the collision operator is self-adjoint with respect to the weight $1/g^\ast$.
\begin{lemma}
The collision operator $\mathcal{L}$, defined in~\eqref{eqn:collision}, is self-adjoint with weight $1/g^\ast$. In particular
\[
\langle\mathcal{L}g\,,h/g^\ast\rangle =  \langle\mathcal{L}h\,,g/g^\ast\rangle\,.
\]
\end{lemma}
\begin{proof}
The proof amounts to direct calculation. Expand the left hand side we have
\[
\begin{aligned}
\langle\mathcal{L}g\,,h/g^\ast\rangle& =\langle \left(-\omega g+\omega g^\ast\frac{\langle\omega g\rangle}{\langle\omega g^\ast\rangle}\right)h/g^\ast\rangle = -\langle\omega gh/g^\ast\rangle + \frac{\langle\omega g\rangle\langle\omega h\rangle}{\langle\omega g^\ast\rangle}\\
& =\langle \left(-\omega h+\omega g^\ast\frac{\langle\omega h\rangle}{\langle\omega g^\ast\rangle}\right)g/g^\ast\rangle = \langle\mathcal{L}h\,,g/g^\ast\rangle
\end{aligned}
\]
which concludes the lemma.
\end{proof}

Now we make the Fr\'echet derivative explicit in the following theorem.

\begin{theorem}\label{thm:frechet}
For a fixed $\gamma = (i,j)$, the Fr\'echet derivative $\nabla_\eta \mathcal{M}_\gamma$ at $\eta$ can be computed as:
\begin{equation}\label{eqn:thm}
\nabla_\eta\mathcal{M}_\gamma=\frac{1}{\langle \omega g^*\rangle}\int_{\mu<0,t}\mu\omega h(x=1,\mu,\omega,t)g_{0}(x=1,-\mu,\omega,t)/g^*d\mu dt
\end{equation}
where $g_{0}$ is the solution to~\eqref{eqn:linear_PTE} with refection coefficient $\eta$ and boundary condition ${\phi_i}$, and $h$ is the solution to the adjoint equation with $\psi_j$ as the boundary condition:
\begin{equation}\label{eqn:adjoint_linear_PTE}
\begin{cases}
\partial_th+\mu \omega\partial_xh&=-\mathcal{L}{h}\\
h(x,\mu,\omega,t=t_{\max})&=0\\
h(x=0,\mu,\cdot)&=\frac{\psi_j(t)g^\ast}{\mu}\,,\mu<0\\
h(x=1,\mu,\cdot)&=\eta h(x=1,-\mu,;)\,,\mu>0\
\end{cases}\,.
\end{equation}
\end{theorem}

\begin{proof}
This is a standard calculus of variation argument. Let $\eta$ be perturbed by $\delta \eta$. Denote the corresponding solution to~\eqref{eqn:linear_PTE} with incoming data $\phi_{i}$ by $g$, and let $\tilde{g} = g-g_0$, then since $g$ and $g_0$ solve ~\eqref{eqn:linear_PTE} with the same incoming data but different reflection coefficients, one derives the equation for $\tilde{g}$:
\begin{equation}\label{eqn:linear_PTE22}
\begin{cases}
\partial_t \tilde{g}+\mu \omega\partial_x \tilde{g}=\mathcal{L}{\tilde{g}}\\
\tilde{g}(x,\mu,\omega,t=0)=0\\
\tilde{g}(x=0,\mu,\cdot)=0\,,\quad \mu>0\\
{\tilde{g}(x=1,\mu,\cdot)=\eta\tilde{g}(x=1,-\mu,\cdot)+\delta\eta g_0(x=1,-\mu,\cdot)}\,,\quad \mu<0
\end{cases}\,,
\end{equation}
where we ignored the higher order term $\delta\eta\delta g$. In this equation, $g_{0}$ serves as the source term at the boundary.

{Now we multiply $\tilde{g}$ equation~\eqref{eqn:linear_PTE22} with $h/g^\ast$ and multiply $h$ equation~\eqref{eqn:adjoint_linear_PTE} with $\tilde{g}/g^\ast$, integrate with respect to $\mu$ and $\omega$ and add them up.} Realizing the $\mathcal{L}$ operator is self adjoint with respect to $1/g^\ast$, the right hand side cancels out and we obtain:
\[
\partial_t\langle \tilde{g}h/g^\ast\rangle+\partial_x\langle \mu\omega \tilde{g}h/g^\ast\rangle= 0\,.
\]
Integrate further in time and space of this equation. Noticing that $\tilde{g}$ has trivial initial data and $h$ has trivial final state data, the $\partial_t$ term drops out and we have:
\[
\int_{t}\langle \mu\omega \tilde{g}h/g^\ast\rangle\rd{t} (x=0)= \int_{t}\langle \mu\omega \tilde{g}h/g^\ast\rangle\rd{t} (x=1)\,.
\]
At $x=0$ plugging in $h(t,x,\mu<0,\omega)=\frac{\psi_j g^\ast}{\mu}$, we have
\[
\begin{aligned}
\int_{t}\langle \mu\omega \tilde{g}h/g^\ast\rangle\rd{t}|_{x=0} &= \int_{t=0}^{t_\text{max}}\int_{\omega=0}^{\omega_\text{max}}\int_{\mu=-1}^0\mu\omega{\frac{1}{g^\ast}\tilde{g}\frac{\psi_jg^\ast}{\mu}}\rd\mu\rd\omega\rd t|_{x=0}\\
& = \int_{t=0}^{t_\text{max}}\psi_j(t)\langle\omega \tilde{g}\rangle\rd t|_{x=0}\\
& = \int_{t=0}^{t_\text{max}}\psi_j(t)\langle\omega g\rangle\ \rd t|_{x=0}- \int_{t=0}^{t_\text{max}}\psi_j(t)\langle\omega g_0\rangle \rd t|_{x=0}\\
&= \left(\mathcal{M}_{ij}(\eta+\delta\eta)-\mathcal{M}_{ij}(\eta)\right)\langle\omega g^\ast\rangle
\end{aligned}
\]
where we used the definition of $\Delta T$.

At $x=1$ we note that the terms with $\eta$ are all canceled out using the reflection condition and the integral reduces to
\[
\int_{t}\langle \mu\omega \tilde{g}h/g^\ast\rangle\rd{t}|_{x=1}=\int_\omega\delta\eta\int_t\int_{\mu=-1}^0\mu\omega h g_0(-\mu)/g^\ast\rd\mu\rd{t}\rd\omega|_{x=1}\,.
\]
Let $\delta\eta\to0$, the two formula lead to:
\[
\nabla_\eta\mathcal{M}_{ij} = \frac{1}{\langle\omega g^\ast\rangle}\int_{t,\mu<0}\mu\omega h(t,x=1,\mu,\omega) g_0(t,x=1,-\mu,\omega)/g^\ast(\omega)\rd\mu\rd{t}\,.
\]
\end{proof}

\begin{remark}
We note that the derivation is completely formal. Indeed the regularity of the solution is not yet known when the incoming data for $h$ has a singularity at $\mu = 0$, making it not even $L_1$. Numerically one can smooth out the boundary condition by adding a mollifier. That is to convolve the boundary condition with a narrow Gaussian term. Numerically, this direct delta will be replaced by a kronecker delta.
\end{remark}

We now summarize the SGD algorithm in Algorithm~\ref{alg:linear}.

 
\begin{algorithm}
\caption{SGD applied on the minimization problem~\eqref{eqn:optimization}.}
\label{alg:linear}
\begin{algorithmic}
\STATE{\textbf{Data:}$IJ$ experiments with \begin{itemize}\item[1.] incoming data $\phi_i$ for $\{i=1\,,\cdots,I\}$;\item[2.] outgoing measurements $\psi_j$ for $\{j=1\,,\cdots,J\}$;\item[3.] error tolerance $\varepsilon$;\item[4.] initial guess $\eta_0$.\end{itemize}}
\STATE{\textbf{Outcome:}The minimizer $\eta_\ast$ to the optimization problem~\eqref{eqn:optimization} that is within $\varepsilon$ accuracy.}
\WHILE{$|\eta_{n+1}-\eta_n|>\varepsilon$}
\STATE{$n=n+1$}
\STATE{Step I: randomly uniformly pick $\gamma_n=(i_n,j_n)\in[1:I]\times[1:J]$;}
\STATE{Step II: compute the solution to the forward problem~\eqref{eqn:linear_PTE} using $\phi_{i_n}$ and $\eta_n$;}
\STATE{Step III: compute the solution to the adjoint problem~\eqref{eqn:adjoint_linear_PTE} using $\psi_{j_n}$ and $\eta_n$;}
\STATE{Step IV: compute $\nabla_\eta L_{\gamma_n}$ using ~\eqref{eqn:thm};}
\STATE{Step V: compute $L_{\gamma_n}$ according to its definition~\eqref{eqn:def_L};}
\STATE{Step VI: update $\eta$ using~\eqref{eqn:SGD} with a step size $\alpha_n$;}
\ENDWHILE
\RETURN $\eta_{n+1}$
\end{algorithmic}
\end{algorithm}

\section{Properties of the equation and the minimization process}\label{sec:properties}
In this section we study some properties of the equation, and the convergence result of the associated optimization problem using SGD. In particular, since the equation is a typical kinetic model with a slightly modified transport term and a linear BGK type collision operator, one would expect some good properties that hold true for general kinetic equations are still valid in this specific situation. Below we study the maximum principle first, and this gives us a natural $L_\infty$ bound of the solution. The result will be constantly used in showing SGD convergence.

\subsection{Maximum Principle}\label{sec:property_MP}
Maximum principle is the property that states the solution in the interior is point-wise controlled by the size of boundary and initial data. It is a classical result for diffusion type equations. In the kinetic framework, it holds true for the radiative transfer equation that has a linear BGK collision operator. But it is not true for the general BGK type Boltzmann equation. If the initial/boundary data is specifically chosen, the Maxwellian term can drag the solution to peak with a value beyond $L_\infty$ norm of the given data. In these cases, an extra constant is expected.

We assume the initial condition to be zero and rewrite the equation as:

\begin{equation}\label{eqn:linear_PTE_steady}
    \begin{cases}
    &\partial_tg+\mu \omega \partial_xg=-\omega g+\frac{\omega g^\ast}{\langle \omega g^\ast\rangle} \langle \omega g\rangle\\
    & g(x=0,\mu,\omega,t) = \phi(\mu,\omega,t)\,,\quad \mu>0\\
    & g(x=1,\mu,\omega,t) = \eta(\omega)g(x=1,-\mu,\omega)\,,\quad \mu<0
    \end{cases}\,,
\end{equation}
where $\eta(\omega)$ is a function bounded between $0$ and $1$. We arrive at the following theorem.

\begin{theorem}\label{thm:maximum_principle}
The equation~\eqref{eqn:linear_PTE_steady} satisfies the maximum principle, namely, there exists a constant $C$ independent of $\phi$ so that
\[
\|g\|_\infty\leq C\|\phi\|_\infty\,.
\]
\end{theorem}
\begin{proof}
The proof of the theorem comes from direct calculation. We first decompose the equation~\eqref{eqn:linear_PTE_steady}. Let $g=g_1+g_2$ with $g_1$ being the ballistic part with the boundary condition and $g_2$ being the remainder, one can then write:
\begin{equation}\label{eqn:linear_PTE_steady_1}
    \begin{cases}
    &\partial_tg_1+\mu\omega\partial_x g_1 = - \omega g_1\\
    & g_1(x=0,\mu>0,\omega,t) = \phi(\mu,\omega,t)\\
    & g_1(x=1,\mu<0,\omega,t) = \eta(\omega)g_1(x=1,-\mu,\omega,t)
    \end{cases}\,,
\end{equation}
and
\begin{equation}\label{eqn:linear_PTE_steady_2}
    \begin{cases}
    &\partial_tg_2+\mu\omega\partial_x g_2 = - \omega g_2+M\langle \omega g\rangle\\
    & g_2(x=0,\mu>0,\omega,t) = 0\\
    & g_2(x=1,\mu<0,\omega,t) = \eta(\omega)g_2(x=1,-\mu,\omega,t)
    \end{cases}\,.
\end{equation}
Note that in this decomposition, $g_1$ and $g_2$ satisfy the same reflective boundary condition at $x=1$, but $g_1$ serves as a source term in the collision term of $g_2$. Both equations can be explicitly calculated. Define the characteristics:
\[
    \frac{dX}{ds}=\mu\omega\,, \quad\text{and let}\quad g_1(t)=g_1(t,X(t))\,,\quad g_2(t)=g_2(t,X(t))\,,
\]
then for every fixed $\mu$ and $\omega$, along the trajectory we have:

\begin{equation}
    \frac{dg_1}{dt}=-\omega g_1\,,\quad    \frac{dg_2}{dt}=-\omega g_2+\frac{\omega g^\ast}{\langle\omega g^\ast\rangle}\langle \omega g\rangle\,,
\end{equation}
which further gives

\begin{equation}
    g_1= g_1(t_0,X(t_0))e^{-\omega(t-t_0)}
\end{equation}
    and
\begin{equation}
    g_2=g_2(t_0,X(t_0))e^{-\omega(t-t_0)}+M\int_{s=t_0}^t\langle\omega g\rangle(s,X(s))e^{-\omega(t-s)}ds\,.
\end{equation}

For $\mu>0$ the characteristic propagates to the right. Depending on the value for $t$ and $x$, the trajectory, when traced back, could either hit the initial data or the boundary data. Let $X(t) = x$, and consider $g=g_1+g_2$, we have:
\begin{enumerate}
    \item If the trajectory hits the boundary data, for $t_0=t-\frac{x}{\mu\omega}>0$ we have $X(t_0)=0$, and
    \begin{equation}\label{eqn:boundary_hit}
        g(t,x)=\phi(t_0)e^{-\omega(t-t_0)}+M\int_{s=t_0}^t\langle\omega  g\rangle(s,X(s))e^{-\omega(t-s)}ds\,.
    \end{equation}
\item If the trajectory hits the initial data, then $t_0=0$ and $X(0)>0$, with the solution being
\begin{equation}\label{eqn:initial_hit}
        g(t,x)=M\int_{s=0}^t\langle\omega  g\rangle(s,X(s))e^{-\omega(t-s)}ds\,.
    \end{equation}
\end{enumerate}
Similarly, for $\mu<0$, if the trajectory hits the boundary data, meaning for $t_1=t-\frac{1-x}{-\mu\omega}>0$ we have $x(t_1)=1$. Define $t_2=t-\frac{2-x}{-\mu\omega}$:

\begin{equation}\label{eqn:g_boudnary_control2}
g=\eta(\omega)\left[\phi(t_2)e^{-\omega(t-t_2)}+Me^{-\omega t}\int_{s=t_2}^{t_1}\langle\omega g\rangle e^{\omega s}ds\right]+Me^{-\omega t}\int_{s=t_1}^{t}\langle\omega g\rangle e^{\omega s}ds\,.
\end{equation}

We multiply~\eqref{eqn:boundary_hit} with $\omega$ and use $\eta<1$ to have:
\begin{equation}\label{eqn:g_boudnary_control}
  |\omega g|\le \omega\|\phi\|_\infty e^{-\omega(t-t_0)}+M\|\langle \omega g\rangle\|_{\infty}[1-e^{-\omega(t-t_0)}]\,.
\end{equation}
Similar inequality can be derived for~\eqref{eqn:initial_hit} as well. Since the bound of~\eqref{eqn:initial_hit} has zero dependence on $\phi$, we ignore its contribution. Taking the moments on both sides of~\eqref{eqn:g_boudnary_control}, we have:
\begin{equation}\label{eqn:max_p}
\langle \omega g\rangle_\text{pos} \le \|\phi\|_\infty\langle\omega e^{-\omega(t-t_0)}\rangle_\text{pos}+\frac{1}{\langle \omega g^\ast\rangle}\|\langle \omega g\rangle\|_{\infty}\langle \omega g^\ast[1-e^{-\omega(t-t_0)}]\rangle_\text{pos}
\end{equation}
where we used the notation $\langle\rangle_\text{pos} = \int_{\omega=0}^\infty\int_{\mu>0}d\mu d\omega$. Similarly, multiplying ~\eqref{eqn:g_boudnary_control2} with $\omega$ and taking moments gives

\begin{equation}\label{eqn:max_n}
\langle \omega g\rangle_\text{neg} \le \|\phi\|_\infty\langle\omega e^{-\omega(t-t_2)}\rangle_\text{neg}+\frac{1}{\langle \omega g^\ast\rangle}\|\langle \omega g\rangle\|_{\infty}\langle \omega g^\ast[1-e^{-\omega(t-t_2)}]\rangle_\text{neg}\,.
\end{equation}
Here $\langle\rangle_\text{neg} = \int_{\omega=0}^\infty\int_{\mu<0}d\mu d\omega$.

Summing up ~\eqref{eqn:max_p} and ~\eqref{eqn:max_n} we have:

\begin{equation}
\langle\omega g\rangle\leq \frac{C_1\|\phi\|_\infty}{C_2}\,,
\end{equation}
with $C_1= \langle\omega e^{-\omega(t-t_0)}\rangle_\text{pos}+\langle\omega e^{-\omega(t-t_2)}\rangle_\text{neg}$, $C_2=\frac{1}{\langle\omega g^\ast\rangle}[\langle \omega g^\ast e^{-\omega(t-t_0)}\rangle_\text{pos}+\langle \omega g^\ast e^{-\omega(t-t_2)}\rangle_\text{neg}]$.

Plug this back into~\eqref{eqn:boundary_hit} or~\eqref{eqn:g_boudnary_control2} we finally have:
\begin{equation*}
    g\leq \Big(1+\frac{C_1}{C_2}\Big)\|\phi\|_\infty=C\|\phi\|_\infty\,,
\end{equation*}
which concludes the proof.
\end{proof}
\subsection{Convergence for noisy data}\label{sec:property_SGD}
 
The study of convergence of SGD is a fairly popular topic in recent years, especially in the machine learning community~\cite{bottou2018optimization, li2019convergence, shamir2013stochastic}. It has been proved that with properly chosen time stepsizes, convergence is guaranteed for convex problems. For nonconvex problems, the situation is significantly harder, and one only seeks for the points where $\nabla f = 0$ where $f = \frac{1}{n}\sum_i f_i$ is the objective function. However, it is hard to distinguish local and global minimizers, and the point that makes $\nabla f=0$ could also be a saddle point, or even a maximizer. See recent reviews in~\cite{bottou2018optimization}.

The situation is slightly different in the PDE-constrained minimization problems. In this setup, the objective function enjoys a special structure: every $f_i$ is the square of the mismatch between one forward map and the corresponding data point, namely $f_i = |\mathcal{M}_i-d_i|^2$. It is typically unlikely for the forward map $\mathcal{M}_i$ to be convex directly, but the outer-layer function is quadratic and helps in shaping the Hessian structure. In particular, since:
\[
\nabla_\eta f_i =2\left(\mathcal{M}_i-d_i\right)\nabla_\eta \mathcal{M}_i\,,\quad\text{and}\quad H_\eta (f_i)=2\nabla_\eta\mathcal{M}_i\otimes\nabla_\eta\mathcal{M}_i +2\left(\mathcal{M}_i-d_i\right) H_\eta(\mathcal{M}_i)\,,
\]
where $H_\eta$ stands for the Hessian term. Then due to the $\nabla_\eta\mathcal{M}_i\otimes\nabla_\eta\mathcal{M}_i$ term, the $H_\eta (f_i)$ would have a better hope of being positive definite, as compared to $H_\eta(\mathcal{M}_i)$ especially when the data is almost clean.

Such structure changes the convergence result. Indeed, in~\cite{jin2020convergence} the authors investigated the performance of SGD applied on problems with PDE constraints. We cite the theorem below (with notation adjusted to fit our setting).

\begin{theorem}
Let $\mathcal{M} = [\mathcal{M}_1\,,\cdots\,,\mathcal{M}_n]$ with $\mathcal{M}_\gamma$ be a forward map that maps $\mathcal{D}(\mathcal{M}_\gamma)\subset X$ to $Y$ with $X$ and $Y$ being two Hilbert spaces with norm denoted by $\|\cdot\|$. Suppose $\mathcal{D}=\cap_\gamma\mathcal{D}(\mathcal{M}_\gamma)$ is non-empty. {Denote $d^\epsilon$ to be the real data perturbed by a noise that is controlled by $\epsilon$}. The loss function is defined as:
\begin{equation}\label{least_squares}
 {   f(\eta)=\frac{1}{n}\sum_{\gamma=1}^n|\mathcal{M}_\gamma(\eta)-d_\gamma^\epsilon|^2}\,.
\end{equation}

Suppose for every fixed $\gamma$:
\begin{enumerate}
    \item The operator $\mathcal{M}_\gamma : X \rightarrow Y$ is continuous, with a continuous and uniformly bounded Fr\'echet derivative on $\mathcal{D}$.
    \item There exists an $\alpha\in (0,\frac{1}{2})$ such that for any $\eta,\tilde{\eta}\in X$,
    \begin{equation}\label{eqn:boundedness}
        |\mathcal{M}_\gamma(\eta)-\mathcal{M}_\gamma(\tilde{\eta})-(\eta-\tilde{\eta})^\top\nabla\mathcal{M}_\gamma(\tilde{\eta})|\le \alpha|\mathcal{M}_\gamma(\eta)-\mathcal{M}_\gamma(\tilde{\eta})|\,.
    \end{equation}

\end{enumerate}
If the Fr\'echet derivative is uniformly bounded in $\mathcal{D}$, and the existence of $\alpha$ holds 
uniformly true for all $\gamma$, then with properly chosen stepsize, SGD converges.

\end{theorem}

We note that in the original theorem, the assumptions are made on $\mathcal{M} =[\mathcal{M}_1\,,\cdots\,,\mathcal{M}_n]$ viewed as a vector. We state the assumptions component-wise but we do require the boundedness and the existence of $\alpha$ to hold true uniformly in $\gamma$. Note also the theorem views $\eta$ as a vector and hence the notation $\eta^\top$ and $\nabla\mathcal{M}$. When viewed as a function of $\omega$, the term reads $\int (\eta-\tilde{\eta})\partial_\eta\mathcal{M}_{\gamma}(\tilde{\eta})\rd\omega$.

To show that the method works in our setting, we essentially need to justify the two assumptions for a properly chosen domain. This is exactly what we will show in Proposition~\ref{prop:Lip_cont} and
Proposition~\ref{prop:cond_2}. In particular, we will show, for each $\gamma=(i,j)$, the Fr\'echet derivative is Lipschitz continuous with the Lipschitz constant depending on $\|\phi_i\|_\infty$ and $\|\psi_j\|_\infty$. By choosing properly bounded $\{\|\phi_i\|_\infty\,,\|\psi_j\|_\infty\}_{ij}$ the Lipschitz constants are bounded. We thus have the uniform in $\gamma$ uniform in $\mathcal{D}$ boundedness of $\partial_\eta\mathcal{M}_{\gamma}$. The second assumption is much more delicate: it essentially requests the gradient to represent a good approximation to the local change and that the second derivative to be relatively weak, at least in a small region. This is shown in Proposition~\ref{prop:cond_2} where we specify the region of such control. 

Now we prove the two propositions.

\begin{proposition}\label{prop:Lip_cont}
For a fixed $\gamma = (i,j)$, $\partial_\eta\mathcal{M}_{\gamma}$ is Lipschitz continuous, meaning there is an $L_\gamma$ that depends on $\|\phi_i\|_\infty$ and $\|\psi_j\|_\infty$ so that
\[
\|\partial_\eta\mathcal{M}_{\gamma}(\eta_1)-\partial_\eta\mathcal{M}_{\gamma}(\eta_2)\|_2\leq L_\gamma\|\eta_1-\eta_2\|_2\,.
\]
\end{proposition}

\begin{proof}
We omit $\gamma=(i,j)$ in the subscript of $\mathcal{M}$ throughout the proof. Recall the definition of $\partial_\eta\mathcal{M}_{\gamma}$, we have:
\begin{equation}
\begin{aligned}
&\partial_\eta\mathcal{M} (\eta_1) - \partial_\eta\mathcal{M}(\eta_2)\\
&=\frac{1}{\langle\omega g^\ast\rangle}\int\limits_{t,\mu<0}\frac{\mu\omega}{g^\ast} \left[h_1(x=1,\mu,\cdot) g_{1}(x=1,-\mu,\cdot)- h_2(x=1,\mu,\cdot) g_{2}(x=1,-\mu,\cdot)\right]\rd\mu\rd{t}
\end{aligned}
\end{equation}
where $h_1$ and $h_2$ solve~\eqref{eqn:adjoint_linear_PTE}, the adjoint equation  equipped with the same boundary and initial data ($\psi_j$). The reflection coefficients are $\eta_1$ and $\eta_2$ respectively. $g_{1,2}$ solve~\eqref{eqn:linear_PTE}, the forward problem, with $\eta_{1,2}$ being the reflective coefficient (also equipped with the same boundary and initial data $\phi_i$).

Call the differences $g=g_1-g_2$ and $h=h_1-h_2$, we rewrite the formula into:

\begin{equation}\label{L cont_1}
    \begin{aligned}
    & \|\partial_\eta\mathcal{M}(\eta_1)-\partial_\eta\mathcal{M}(\eta_2)\|_2^2\\
&=  \int_\omega\frac{\omega^2}{{ g^*}^2\langle \omega  g^*\rangle^2}\Big[\int_{\mu<0,t}\mu[h_{1}(x=1,\mu,\cdot)g(x=1,-\mu,\cdot)+g_{2}(x=1,-\mu,\cdot)h(x=1,\mu,\cdot)]\Big]^2\\
\leq&  C_1\int_\omega\int_{\mu<0,t}\frac{\omega^2\mu^2}{{ g^*}^2}h^2_{1}(x=1,\mu,\cdot)g^2(x=1,-\mu,\cdot)\\ +&C_1\int_\omega\int_{\mu<0,t}\frac{\omega^2\mu^2}{{ g^*}^2}g^2_{1}(x=1,\mu,\cdot)h^2(x=1,-\mu,\cdot)\,
        \end{aligned}
\end{equation}
Here $C_1=\frac{2t_{\max}}{\langle\omega g^\ast\rangle^2}$. Since the two terms are similar, we only treat the first one as an example below. The same derivation can be done for the second term.

According to the definition of $g$, it follows the forward phonon transport equation~\eqref{eqn:linear_PTE2} with trivial boundary and initial data, but with extra source terms. In particular:
\[
\partial_t g +\mu\omega\partial_x g =-\omega g + M\langle\omega g\rangle\,
\]
with trivial boundary condition at $x=0$, trivial initial condition, and 
\[
\begin{aligned}
g(t,x=1,\mu,\omega) &= \eta_1 g_1(t,x=1,-\mu,\omega) -\eta_2  g_2(t,x=1,-\mu,\omega) \\
& = \eta_1g(t,x=1,-\mu,\omega)+\left(\eta_1-\eta_2\right)g_2(t,x=1,-\mu,\omega)\\
& = \eta_1g(t,x=1,-\mu,\omega) + S\,.
\end{aligned}
\]
Here we denote $S = \left(\eta_1-\eta_2\right)g_2(t,x=1,-\mu,\omega)$, and it can be easily seen that $\|S\|_2\leq C\|\eta_1-\eta_2\|_2$. This $C$ depends on $\|g\|_\infty$ which then relies on $\|\phi_i\|_\infty$. To bound~\eqref{L cont_1} amounts to controlling the two terms using $S$.

We can solve the equation explicitly. In particular, for fixed $\mu,\omega$, we use the method of characteristics by defining:
\[
\frac{\rd X}{\rd s} = \mu\omega\,,\quad X(t) = y
\]
and the solution for $g$ can be written down as
\begin{equation*}
    g(t,y,\cdot)=g(t_0,X(t_0),\cdot)e^{-\omega(t-t_0)}+Me^{-\omega t}\int_{t_0}^t\langle\omega g\rangle(s,X(s))e^{\omega s}ds\,.
\end{equation*}

This formula allows us to trace $g$ back to either the initial and boundaries.

In the case of $\mu>0$, the trajectory is right-going, and since initial condition and the boundary at $x=0$ are zero, the solution becomes

\begin{equation}\label{eqn:pos_mu_g_formula}
    g(t, y,\cdot)=Me^{\omega t}\int_{t_3}^t\langle\omega g\rangle(s,x(s))e^{-\omega s}ds\leq\frac{M}{\omega}\|\langle\omega g\rangle\|_{\infty}[1-e^{-\omega(t-t_3)}]\,,
\end{equation}
where we have $t_3 = t-\frac{y}{\mu\omega}$. This leads to
\begin{equation}\label{eqn:pos_mu_g}
\langle\omega g\rangle_\text{pos}\leq \langle{M}[1-e^{-\omega(t-t_3)}]\rangle_\text{pos} \|\langle\omega g\rangle\|_{\infty}\,.
\end{equation}
Here we used the notation $\langle\cdot\rangle_\text{pos} = \int_{\omega,\mu>0}d\mu d\omega$.

In the case of $\mu<0$, the trajectory is left-propagating. Tracing it backwards the trajectory hits $x=1$. Let $t_2=t-\frac{1-y}{(-\mu)\omega}$ and $t_1 = t_2-\frac{1}{(-\mu)\omega}$, we have $X(t_2) = 1$ and $X(t_1) = 0$. Suppose $t_1\geq 0$, the solution reads:
\begin{equation}
    \begin{aligned}
    g(t,y,\mu,\omega)&=g(t_2,1,\mu,\omega)e^{-\omega(t-t_2)}+Me^{-\omega t}\int_{t_2}^t\langle\omega g\rangle(s,X(s))e^{\omega s}ds\\
    &=[\eta_1 g(t_2,1,-\mu,\omega)+S(-\mu,\omega)]e^{-\omega(t-t_2)}+Me^{-\omega t}\int_{t_2}^t\langle\omega g\rangle(s,X(s))e^{\omega s}ds\\
&=e^{-\omega(t-t_2)}S(-\mu,\omega)+ \eta_1 Me^{-\omega t}\int_{t_1}^{t_2}\langle\omega g\rangle(s,x(s))e^{\omega s}ds \\
&+Me^{-\omega t}\int_{t_2}^t\langle\omega g\rangle(s,X(s))e^{\omega s}ds\,.
\end{aligned}
\end{equation}

Noting $\eta_1<1$, we pull out $\|\langle\omega g\rangle\|_\infty$ and have:
\begin{equation*}
\begin{aligned}
   g(t,y,\mu,\omega)&\le e^{-\omega(t-t_2)}S(-\mu,\omega)+M\|\langle\omega g\rangle\|_{\infty}e^{-\omega t}\big[\int_{t_1}^{t_2}e^{\omega s} ds+\int_{t_2}^t e^{\omega s}ds\big]\\
&= e^{-\omega(t-t_2)}S(-\mu,\omega)+\frac{M}{\omega}\|\langle\omega g\rangle\|_{\infty}[1-e^{-\omega(t-t_1)}]\,.
\end{aligned}
\end{equation*}

Note that if $t_i<0$, the lower bound of the integral is replaced by $0$ and the inequalities still hold true. 
\begin{equation}\label{eqn:neg_mu_g}
\langle\omega g\rangle_\text{neg}\leq \langle\omega e^{-\omega(t-t_2)}S(-\mu,\omega)\rangle_\text{neg} +\langle{M}[1-e^{-\omega(t-t_1)}]\rangle_\text{neg} \|\langle\omega g\rangle\|_{\infty}\,.
\end{equation}
Here we used the notation $\langle\cdot\rangle_\text{neg} = \int_{\omega,\mu<0}d\mu d\omega$.

Summing up~\eqref{eqn:pos_mu_g} and~\eqref{eqn:neg_mu_g}, use Cauchy-Schwarz, and $\langle M\rangle_{\mu\lessgtr0,\omega} = \frac{1}{2}\langle M\rangle$, we have:
\begin{equation}\label{eqn:sum}
\langle\omega g\rangle\leq \|S\|_2\|\omega e^{-\omega(t-t_2)}\|_2 +\|\langle\omega g\rangle\|_\infty\left(1- \langle Me^{-\omega(t-t_1)}\rangle_\text{neg}- \langle Me^{-\omega(t-t_3)}\rangle_\text{pos}\right)\,.
\end{equation}

According to the definition of $t_3$, $t-t_3 = \frac{y}{\mu\omega}$, and by setting $t_3>0$, we have $\frac{y}{t\omega}<\mu<1$:
\begin{equation}
    \langle M e^{-\omega(t-t_3)}\rangle_\text{pos}=\int_{\omega=y/t}^{\infty}\int_{\mu=y/\omega t}^{1} Me^{-\frac{y}{\mu}}d\mu d\omega\leq \frac{1}{y}\int_{y/t}^\infty (e^{-y}-e^{-\omega t})M d\omega\,.
\end{equation}

 Similarly, we have $t-t_1= \frac{2-y}{-\mu\omega}$, and by setting $t_1>0$, we have $\frac{2-y}{t\omega}<|\mu|<1$:
 \begin{equation}
   \langle M e^{-\omega(t-t_1)}\rangle_\text{neg}=\int_{\omega=(2-y)/t}^{\infty}\int_{\mu=(2-y)/\omega t}^{1} Me^{-\frac{2-y}{\mu}}d\mu d\omega\,.
\end{equation}

Denote by $\alpha/2$ the minimum of the two equations. It is less than $1$ and strictly positive for $y\in[0,1]$. Plug it back in~\eqref{eqn:sum}, we have
\[
\alpha\|\langle\omega g\rangle\|_\infty\leq \|S\|_2\|\omega e^{-\omega(t-t_2)}\|_2\,.
\]
This estimate, when plugged back in~\eqref{eqn:pos_mu_g_formula} and insert into the first term in~\eqref{L cont_1}, we have, $y = 1$:
\[
\begin{aligned}
&\int_{\omega,\mu<0,t}\frac{\omega^2\mu^2}{(g^\ast)^2}h^2_1(t,y,\mu,\omega)g^2(t,y,-\mu,\omega)d\omega d\mu dt\\
\leq &\int_{\omega,\mu<0,t}\frac{M^2}{(g^\ast)^2}h^2_1(t,y,\mu,\omega)\left(1-e^{-\frac{1}{|\mu|}}\right)^2d\omega d\mu dt\|\langle\omega g\rangle\|_\infty^2\\
\leq &\frac{1}{\alpha^2}\int_{\omega,\mu<0,t}\frac{\omega^2}{\langle\omega g^\ast\rangle^2}h^2_1\left(1-e^{-\frac{1}{|\mu|}}\right)^2d\omega d\mu dt\|S\|_2^2\|\omega e^{-\omega(t-t_2)}\|_2^2\\
\leq & C\|S\|_2^2\leq C\|\eta_1-\eta_2\|_2^2\,.
\end{aligned}
\]
where $C$ depends on $\|\langle\omega h_1\rangle\|_\infty$ which relies on $\|\psi_j\|_\infty$, according to Theorem~\ref{thm:maximum_principle}. This concludes the proof by setting $L_\gamma^2 = 2CC_1$. It depends on $\gamma$ but is independent of the particular choice of $\eta_{1,2}$.
\end{proof}

It is then easy to conclude that by setting $\|\phi_i\|_\infty$ and $\|\psi_j\|_\infty$ to be upper and lower bounded, we have the Lipschitz constant uniformly upper bounded.

\begin{proposition}\label{prop:cond_2}
For a fixed $\gamma =(i,j)$, denote $L_\gamma$ the Lipschitz constant of $\partial_\eta\mathcal{M}_{\gamma}$. Suppose there is a constant $M_\gamma$ such that $|\mathcal{M}_\gamma(\eta_1)-\mathcal{M}_\gamma(\eta_2)|>M_\gamma\|\eta_1-\eta_2\|_2$, then it holds in a small neighborhood of $\eta_\ast$, the optimal solution to~\eqref{eqn:optimization}, with radius $\frac{M_\gamma}{4L_\gamma}$ that:
\begin{equation}\label{eqn:condtion_2_prop}
    |\mathcal{M}_\gamma(\eta_1)-\mathcal{M}_\gamma(\eta_2)-\partial_\eta\mathcal{M}_{\gamma}(\eta_2)^\top(\eta_1-\eta_2)|\le \alpha|\mathcal{M}_\gamma(\eta_1)-\mathcal{M}_\gamma(\eta_2)|\,.
\end{equation}
\end{proposition}
\begin{proof}
We omit the subscript $\gamma$ in the proof. Without loss of generality, let $\mathcal{M}(\eta_1)>\mathcal{M}(\eta_2)\ge0$. To show the proposition amounts to finding an $\alpha\in(0,\frac{1}{2})$ such that

\begin{equation*}
-\alpha[\mathcal{M}(\eta_1)-\mathcal{M}(\eta_2)]\le\mathcal{M}(\eta_1)-\mathcal{M}(\eta_2)-\partial_\eta\mathcal{M}(\eta_2)^\top(\eta_1-\eta_2)\le \alpha[\mathcal{M}(\eta_1)-\mathcal{M}(\eta_2)]\,.
\end{equation*}

The two sides are symmetric and we only prove the second inequality. Noting that from the mean value theorem, there is $p\in(0,1)$ such that
\begin{equation*}
     \mathcal{M}(\eta_1)-\mathcal{M}(\eta_2)=\partial_\eta\mathcal{M}(\eta_2+p(\eta_1-\eta_2))^\top(\eta_1-\eta_2)>0\,,
\end{equation*}
then the inequality translates to
\begin{equation}\label{eqn:prop_bounded}
[\partial_\eta\mathcal{M}(\eta_2+p(\eta_1-\eta_2))-\partial_\eta\mathcal{M}(\eta_2)]^\top(\eta_1-\eta_2)\le\alpha M\|\eta_1-\eta_2\|\,.
\end{equation}
Due to the Lipschitz condition on the Fr\'echet derivaitve:
\[
[\partial_\eta\mathcal{M}(\eta_2+p(\eta_1-\eta_2))-\partial_\eta\mathcal{M}(\eta_2)]^\top(\eta_1-\eta_2)\leq L\|\eta_1-\eta_2\|^2_2\,.
\]
Since $\eta_i$ are chosen from the neighborhood $B(\eta_\ast,\frac{M}{4L})$, then $\|\eta_1-\eta_2\|_2\leq \frac{M}{2L}$ and the inequality holds true with $\alpha = 1/2$.
\end{proof}

\begin{remark}
We do comment that the assumption $|\mathcal{M}_\gamma(\eta_1)-\mathcal{M}_\gamma(\eta_2)|>M_\gamma\|\eta_1-\eta_2\|_2$ is strong. It essentially implies that the Lipschitz constant has a lower bound. If not, then there is a possibility to find $\eta_1\neq\eta_2$ with $\mathcal{M}_\gamma(\eta_1) =\mathcal{M}_\gamma(\eta_2)$. This leads to the unique reconstruction impossible theoretically. As what we emphasized above, theoretically proving the unique reconstruction is beyond the scope of the current paper and we merely state it as an assumption here.
\end{remark}

\section{Numerical Results}\label{sec:numerics}
We present our numerical results in this section. As a numerical setup, we choose $x\in[0,0.5]$, $t\in[0,5]$, $\mu\in [-1,1]$, and $\omega\in[\omega_{\min}, \omega_{\max}]$. Meanwhile, we set the discretization to be uniform with $\Delta x=0.02$, $\Delta t=0.01$, $\Delta \mu=0.01$ and $\Delta \omega=0.05$. In space we use the upwind method, and in $\omega$ and $\mu$ direction we use the discrete ordinate method.

\subsection{Examples of the forward equation}
\begin{figure}[htb]
\centering
\includegraphics[scale=1]{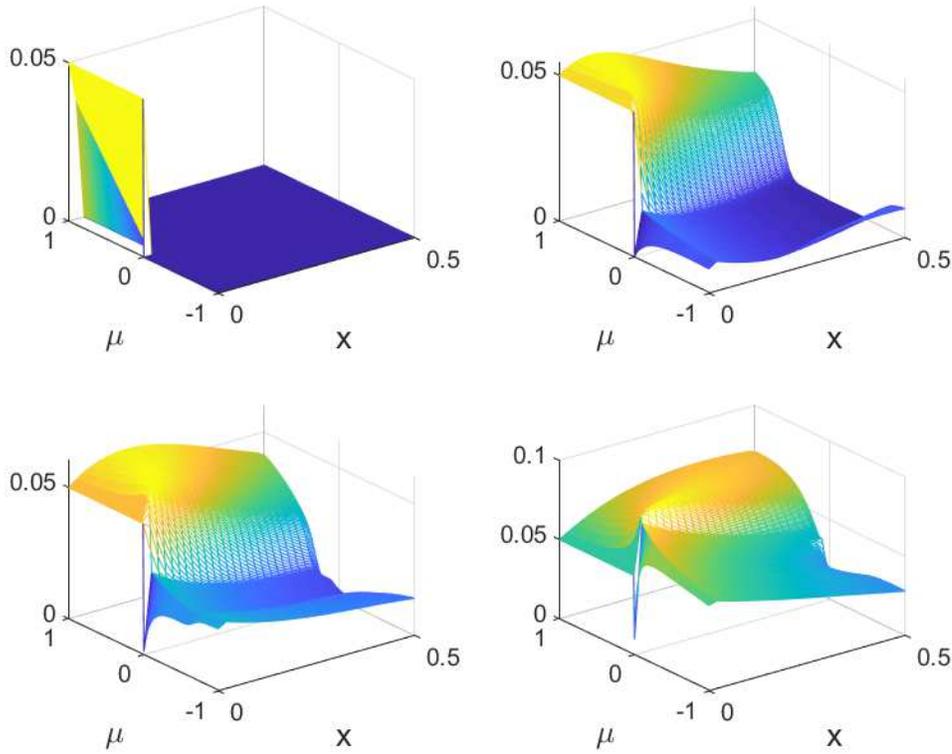} 
\caption{Solution to the forward PDE ($\int gd\omega$ as a function of $(x,\mu)$ at different time frames) with input~\eqref{eqn:boundary_num_forward}. The four plots from top left to bottom right are the solution at time being $0.01$, $0.5$, $1$ and $3$ respectively.}
\label{fig:forward}
\end{figure}

\begin{figure}[htb]
\centering
\includegraphics[scale = 0.75]{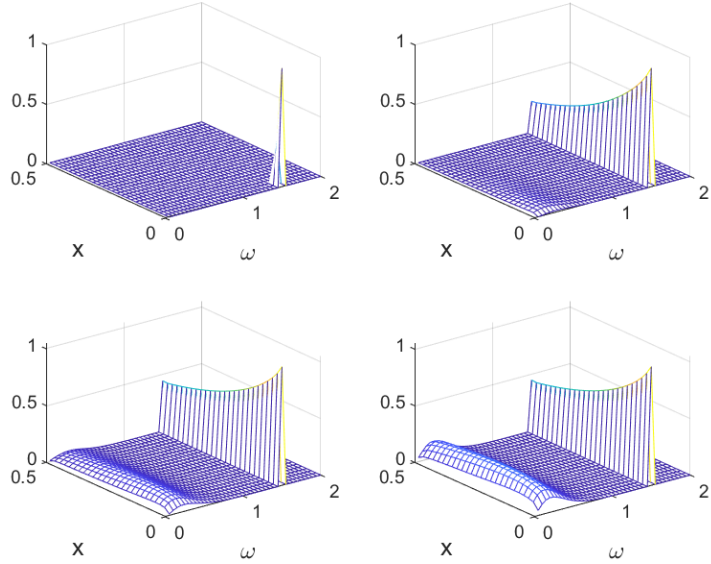} 
\caption{Solution to the forward PDE ($\int gd\mu$ as a function of $(x,\omega)$ at $t=0.01, 0.5, 1.5, 3$ respectively). The boundary input is centered at $\omega=1.5$.}
\label{fig:forwardmu}
\end{figure}

We show an example of the forward solution in Figure~\ref{fig:forward} and Figure~\ref{fig:forwardmu}. {The forward solver is an extended version of an asymptotic preserving algorithm developed in~\cite{li2017implicit}.} In the example, the initial condition is set to be zero, and the boundary condition is set to be 
\begin{equation}\label{eqn:boundary_num_forward}
g(x=0,\mu>0,\omega,t)=\delta(\omega-1.5)\,.
\end{equation}
{Here we implement $\delta$ function as the Kronecker delta function concentrated at the grid point $\omega=1.5$ where it takes the value 1.} As is visualized in Figure~\ref{fig:forward}, the wave gradually propagates inside the domain, with the larger $\mu$ moving with a faster speed. At around $t=0.5$, the wave is reflected back and starts propagating backwards to $x=0$. Besides the transport term $\mu\omega\partial_x$, the BGK collision term ``smooths" out the solution. This explains the nontrivial value to negative $\mu$ even before the reflection taking place. In Figure~\ref{fig:forwardmu} we integrate the solution in $\mu$ and present it as a function on $(x,\omega)$. The solution still preserves the peak at $\omega = 1.5$ as imposed in the boundary condition, but some mass is transported to the smaller $\omega$ values due to the BGK term.

\begin{figure}[htb]
\centering
\includegraphics{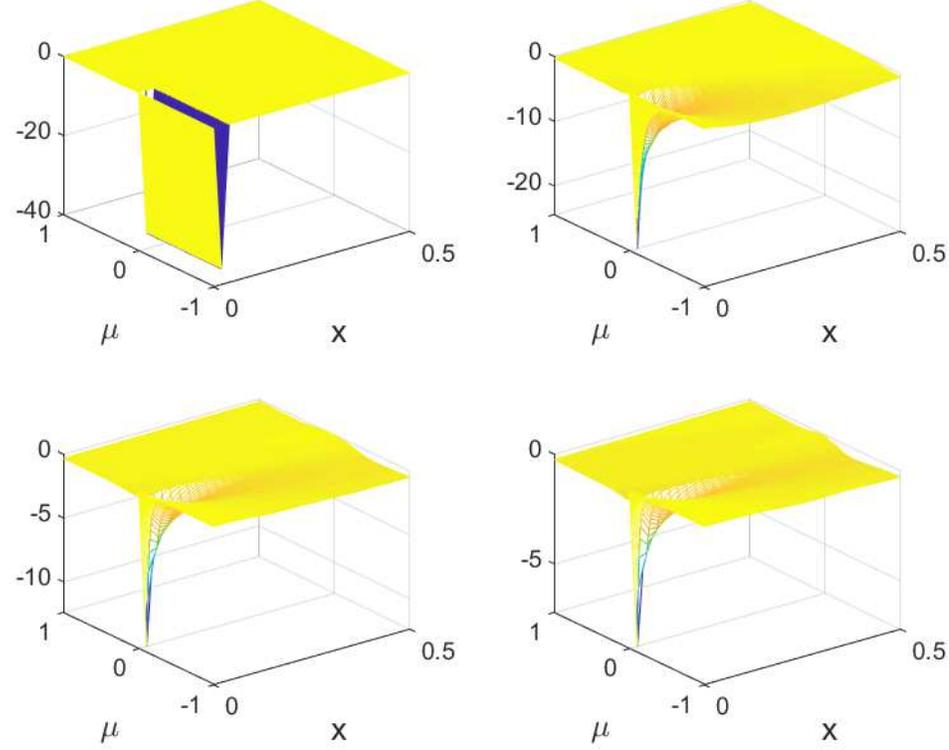} 
\caption{Evolution of the solution to the adjoint PDE. We plot $\int h d\omega$ as a function of $\mu$ and $x$ backwards in time. The time frame for the four plots are $t=5$, $4$, $2$ and $0.01$ respectively.}\label{fig:adjoint}
\end{figure}

Similarly for a demonstration, we also plot an example of solution to the adjoint equation in Figure~\ref{fig:adjoint}. The final and boundary conditions are set to be
\begin{equation}\label{eqn:boundary_num_adjoint}
\begin{aligned}
h(x,\mu,\omega,t\geq t_\text{max})=0\,,\quad & {h(x=1,\mu<0,\omega,t)=M(\omega)\frac{\delta(t-t_{\max})}{\mu\omega\Delta t},}
\end{aligned}
\end{equation}

with $M(\omega)=\frac{(10\omega)^3e^{10\omega}}{(e^{10\omega}-1)^2}$. We note that the plot needs to be viewed backwards in time. The final time data is specified, meaning we need to make sure the solution at $t=t_{\max}=5$ is trivial.

\subsection{Inverse example I}
\begin{figure}[htb]
\centering
\includegraphics[scale = 0.5]{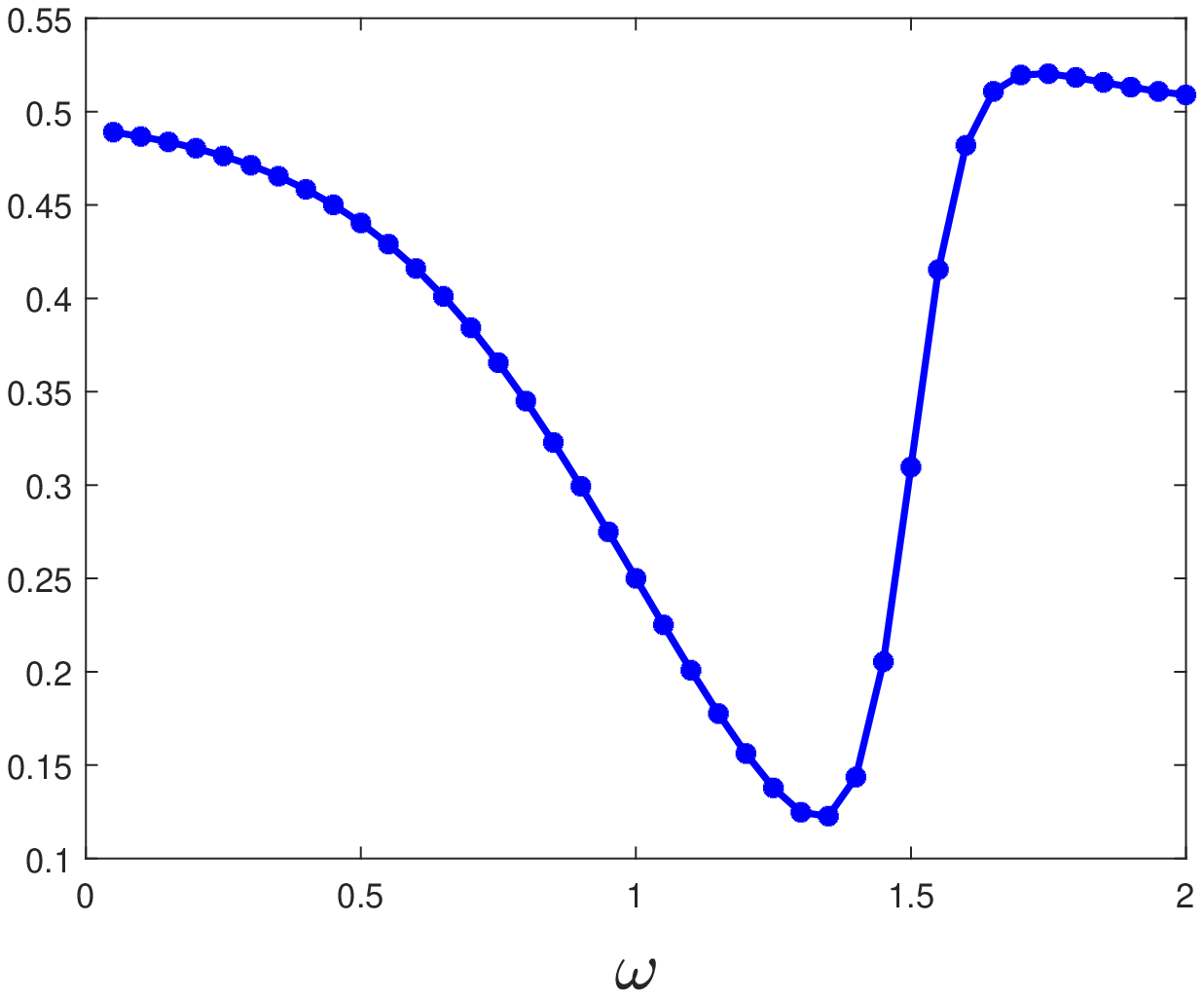} 
\caption{Profile of reference reflection coefficient $\eta_\text{ref}(\omega)$.}\label{fig:ref_eta}
\end{figure}

We now show our solution for the inverse problem. In the first example, we parametrize the reflection coefficient, namely we assume $\eta$ has the form of:
\begin{equation}\label{eqn:num_ref_eta}
\eta=(\tanh(10(\omega-a))-\tanh(2(\omega-b)))/4+1/2\,,
\end{equation}
where $a$ and $b$ are two parameters to be found. We set the ground-truth configuration to be $(a^\ast,b^\ast) = (1.5\,,1)$ and we plot the reference reflection coefficient $\eta_\text{ref}$ in Figure~\ref{fig:ref_eta}. This reference solution is chosen to reflect the investigation in~\cite{hua2017experimental}.

For testing, we set $I=40$ and $J=1$ and use the following $\phi_i$ and $\psi_j$:
\begin{align}\label{eqn:num_phi_psi}
\phi_i(\mu,\omega,t) &= \delta(\omega-\omega_i)=\phi_i(\omega)\,,\quad i=1,\cdots,40\,,\\
\psi_j(\mu,\omega,t)&=\delta(t-t_{\max}), \quad j=J
\end{align}
where $\omega_i = \omega_{\min}+i\Delta\omega$ is the discrete point. {We set $[\omega_{\min},\omega_{\max}]=[0.05,2]$. }

\begin{figure}[htb]
\centering
\includegraphics[scale = 0.75]{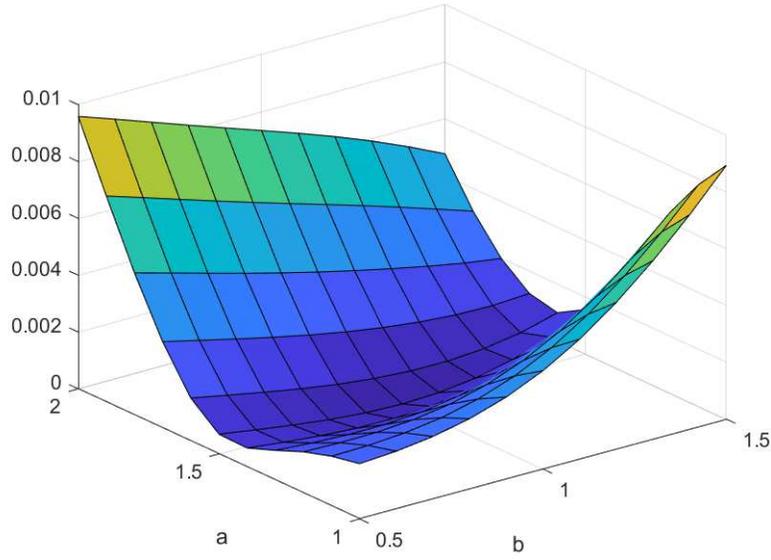} 
\caption{Cost function computed in the neighborhood of $(a^\ast,b^\ast)$. We choose $(a,b)\in [1,2]\times[0.5,1.5]$ with uniform mesh $0.1$. It is clear the cost function is convex in $(a,b)$.}\label{fig:costfunction}
\end{figure}

In Figure~\ref{fig:costfunction} we plot the cost function $L$ in the neighborhood of $(a^\ast,b^\ast)$. {It can be seen that the cost function is convex in both parameters (see Appendix A for the proof). We note that this is the case for this special example where $\eta$ is parameterized. In reality, however, one typically does not have the profile of the ground-truth in hand to make an accurate prediction of the parameterization form.}

In this parametrized setting, the SGD algorithm is translated to

\begin{equation}
a_{n+1}=a_n-2\alpha L_{i_nj_n} G_{i_nj_n,a}\,,\quad b_{n+1}=b_n-2\alpha L_{i_nj_n} G_{i_nj_n,b}\,,
\end{equation}
where $(i_n,j_n)$ are randomly drawn from $[1:I]\times [1,J]$. Here the Fr\'echet derivative can be explicitly computed:
\begin{equation}
\begin{aligned}
G_{ij,a}&= \frac{\partial L}{\partial\eta}\frac{\partial\eta}{\partial a}=-\int_{\omega}10(1-\tanh^2(10(\omega-a_n)))G_{ij}(\omega)/4d\omega\,,\\
G_{ij,b}&= \frac{\partial L}{\partial\eta}\frac{\partial\eta}{\partial b}=\int_{\omega}2(1-\tanh^2(2(\omega-b_n)))G_{ij}(\omega)/4d\omega\,.
\end{aligned}
\end{equation}
with $G_{ij}(\omega) = \partial_\eta L_{ij}$ standing for the Fr\'echet derivative of $L_{ij}$ with respect to $\eta$ and can be computed according to Theorem~\ref{eqn:thm}. 

\begin{figure}
\centering
\includegraphics[scale = 0.5]{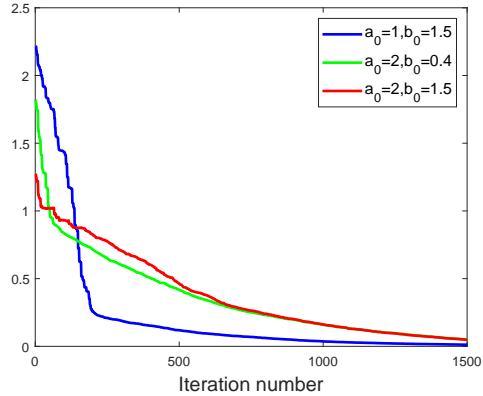}
\caption{The decay of $L_2$ norm of the error in the reconstruction process with three different initial configurations.}\label{fig:decay}
\end{figure}

We run the SGD algorithm with three different initial guesses {with $(a_0,b_0)$ being $(1,1.5)$, $(2,0.4)$ and $(2,1.5)$}. The decay of the error in the reflection coefficient is plotted in Figure~\ref{fig:decay}. Here error is defined to be
\[
{\text{error} = \left(\sum\limits_{i}[\eta_{(a_n, b_n)}(\omega_i)-\eta_{\text{ref}}(\omega_i)]^2\right)^{1/2}}\,.
\]
{where $\eta_{(a_n, b_n)}$ denotes the reflection coefficient for parameters $(a_n,b_n)$.}

\subsection{Inverse Example II}
In this example, we use exactly the same configuration as in the previous example with the same reference solution~\eqref{eqn:num_ref_eta} with $(a^\ast\,,b^\ast)=(1.5\,,1)$. $\phi_i$ and $\psi_j$ are also defined in the same way as in~\eqref{eqn:num_phi_psi}. The main difference here is that we do not parameterize $\eta$: It is viewed as a completely unknown vector of $40$ dimensions (with $\omega_{\max}=2$ and $\omega_{\min}=0.05$ and $\Delta \omega = 0.05$).

\begin{figure}[htb]
\centering
\includegraphics[width =0.32\textwidth, height = 0.2\textheight]{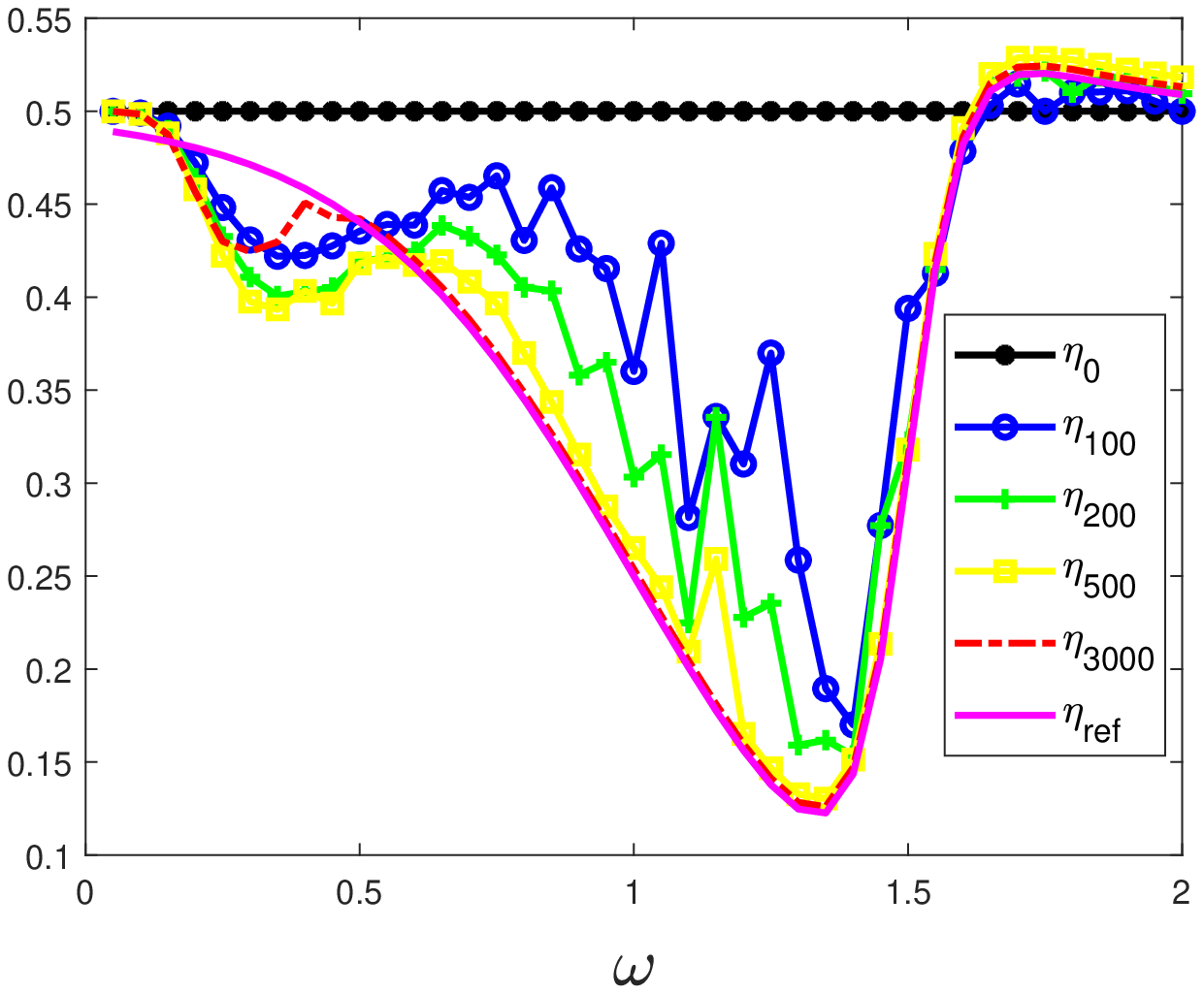}
\includegraphics[width =0.32\textwidth, height = 0.2\textheight]{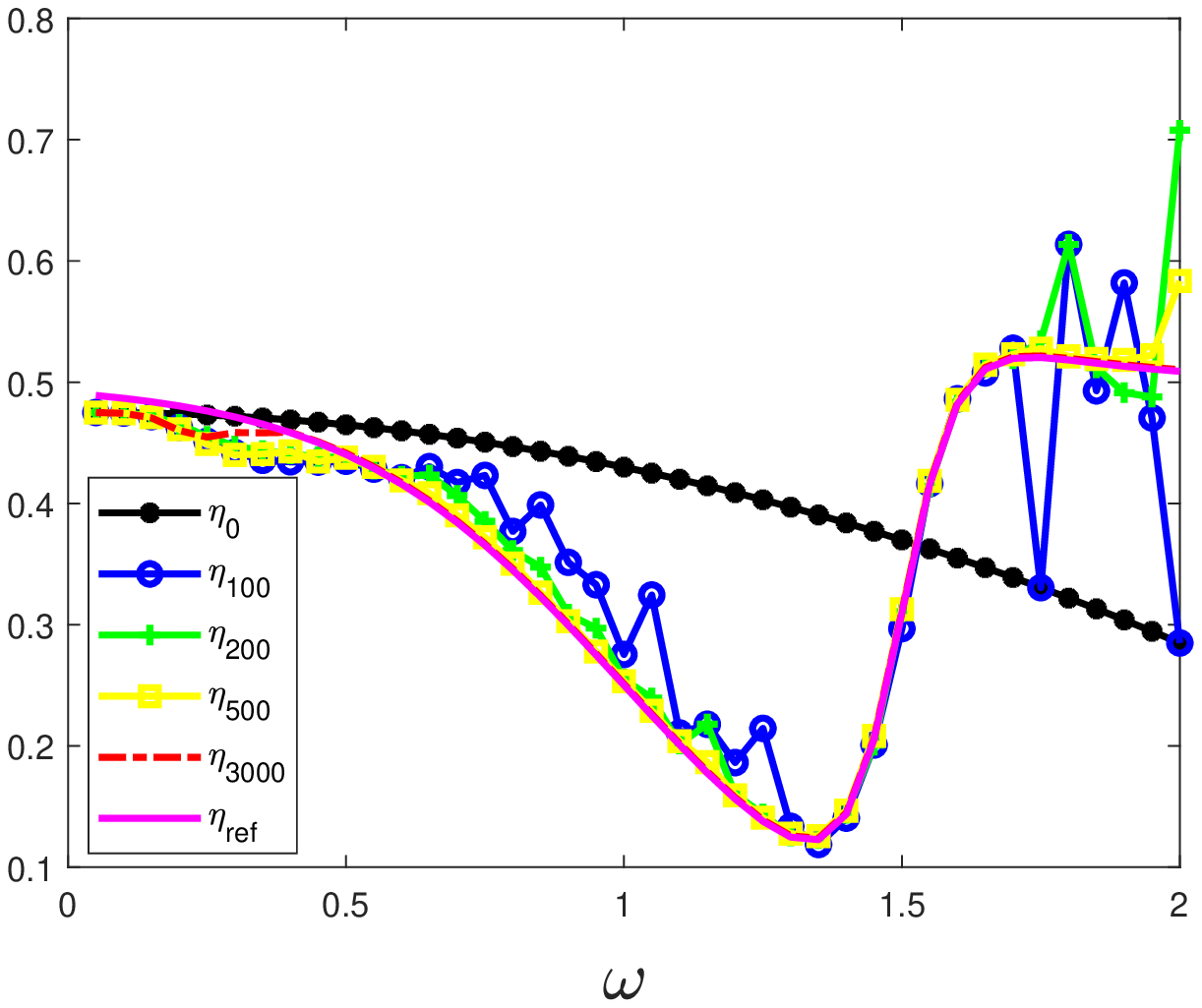}
\includegraphics[width =0.32\textwidth, height = 0.2\textheight]{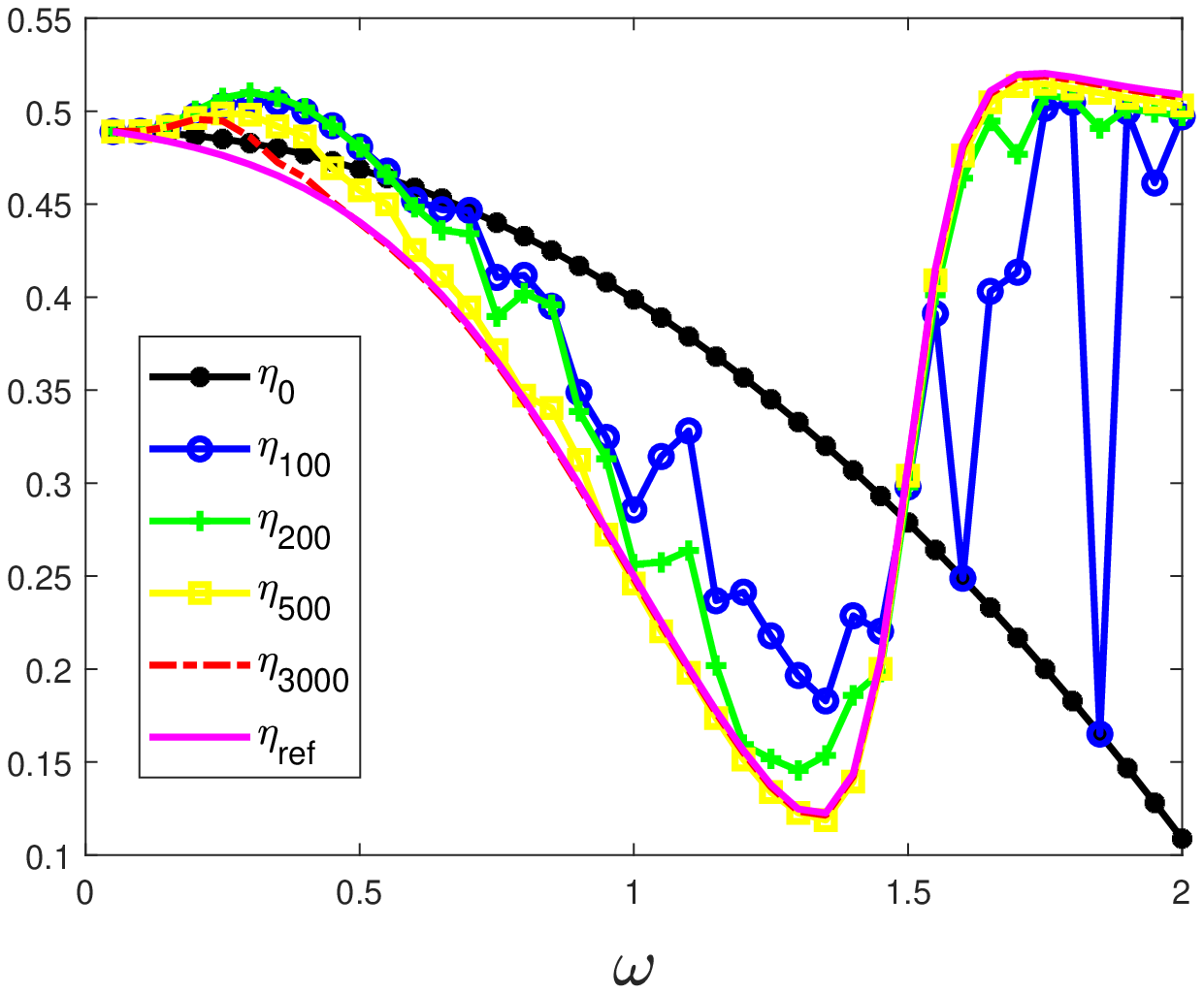}
\caption{Profile of reflection coefficient at different iteration steps. The initial guess for the three plots are $\eta_0=0.5$, $\eta_0=0.4750-0.05(\omega-0.05)^2$ and $\eta_0=0.4891-0.1(\omega-0.05)^2$ respectively. The solution at different iterations are plotted with different colors ($\eta_n(\omega)$ with $n=100$, $200$, $500$ and $3000$ are presented by blue, green, yellow and red lines). The reference solution is shown as the pink line. In all three plots, the solution at $n=3000$ almost recovers the reference solution.}\label{fig:decay1}
\end{figure}

\begin{figure}[h]
\centering
\includegraphics[width=0.5\linewidth]{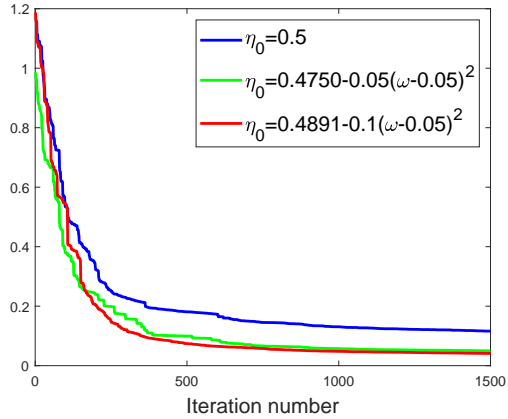}
\caption{The decay of error in terms of the number of iteration for three different initial guesses. The error saturates after about $1500$ iterations.}
\label{fig:error_ex2}
\end{figure}

{We run SGD algorithm with three sets of initial guesses given by $\eta_0=0.5$, $\eta_0=0.4750-0.05(\omega-0.05)^2$ and $\eta_0=0.4891-0.1(\omega-0.05)^2$. The optimization results at different time slices are plotted in Figure~\ref{fig:decay1}, and the decay of $L_2$ norm of error is also plotted in Figure~\ref{fig:error_ex2}. As can be seen in these plots, the numerical solution to the optimization problem reconstructs the ground-truth reflection coefficients. We note that the error in Figure~\ref{fig:error_ex2} exponentially decay at the beginning and eventually saturates. This comes from the numerical error from discretization. In the computation of the Fr\'echet derivation, we unavoidably introduce discretization error. This error cannot be overcome with fixed stepsize in optimization, see~\cite{bottou2018optimization} Theorem 4.6 and 4.8 with non-zero $M$.}

\begin{figure}
\centering
\includegraphics[scale = 0.5]{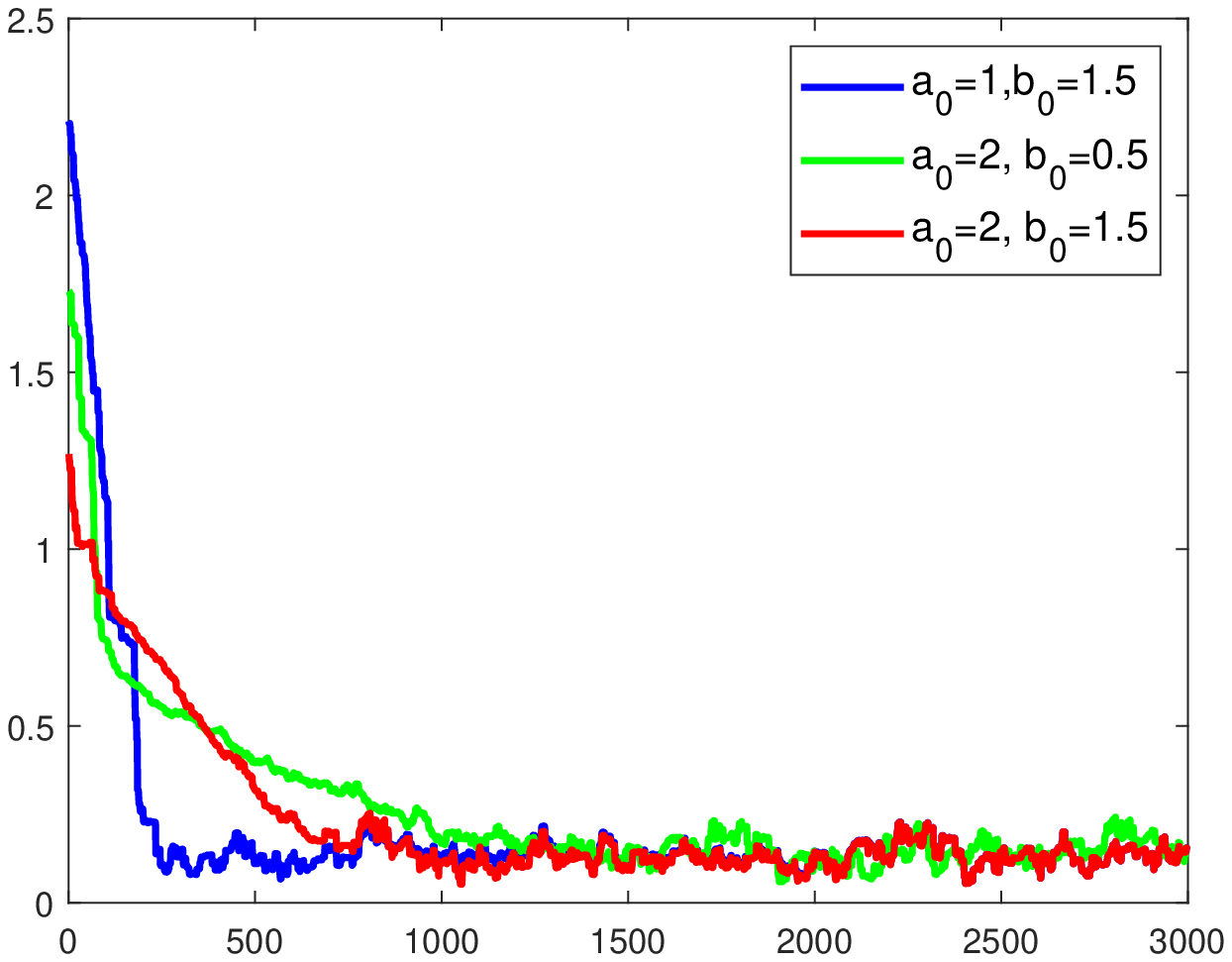}
\caption{{The decay of $L_2$ norm of the error in the reconstruction process with three different initial configurations in the presence of noisy data.}}\label{fig:decay2}
\end{figure}

\begin{figure}
\centering
\includegraphics[width =0.45\textwidth, height = 0.2\textheight]{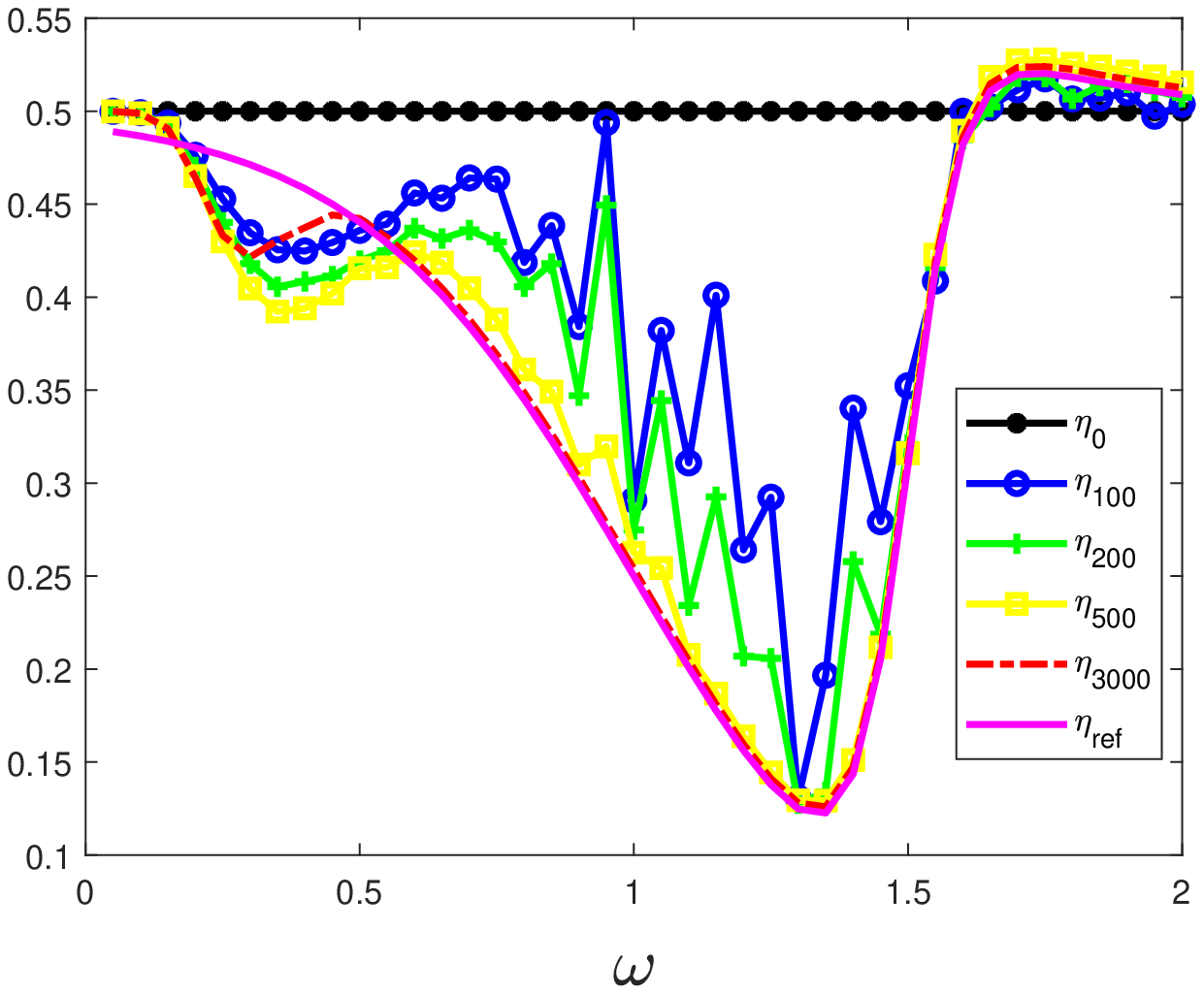}
\includegraphics[width =0.45\textwidth, height = 0.2\textheight]{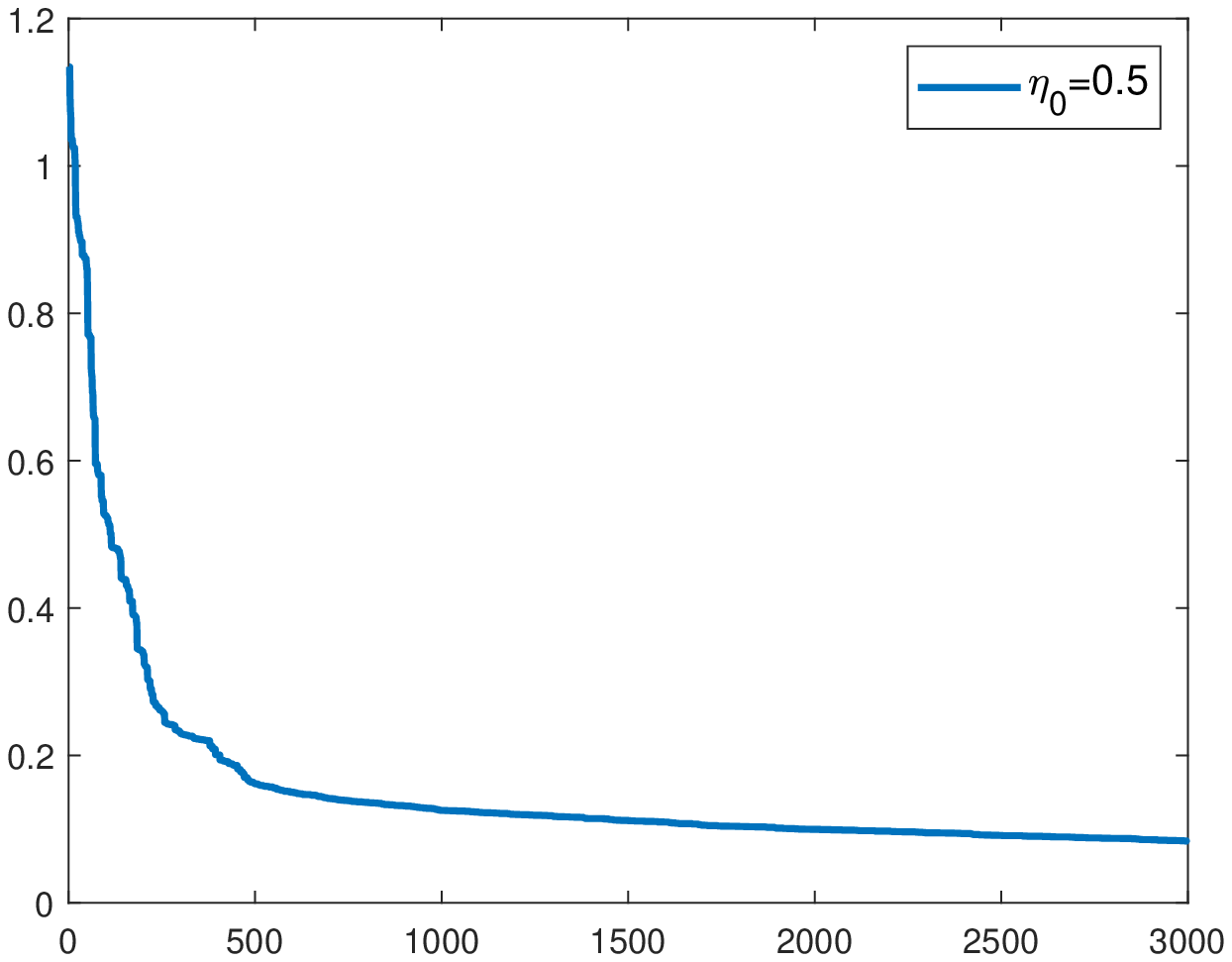}
\caption{{Reconstruction of the reflection coefficient at different iterations when $J=50$. We show the profile of the reconstructed reflection coefficient at different stages of the iteration and the decay of the error.}}\label{fig:decay3}
\end{figure}

{\subsection{Inverse examples III}
In the last example, we repeat the numerical experiment for the previous examples with added noise. The noise is chosen as a random variable uniformly sampled with noise level $2.5\%$.
As in Example I, we run the SGD algorithm with three different initial guesses with $(a_0,b_0)$ given by $(1,1.5)$, $(2,0.5)$ and $(2,1.5)$. The decay of the error in the reflection coefficient is plotted in Figure~\ref{fig:decay2}. In all the three cases, the error saturates after about 1500 iterations.}\\
{
Finally we repeat the numerical experiments with added noise to Example II. In this case we set $J=50$, with measuring operators set as delta functions centered at $t_j$ that are randomly selected in the time interval $[4.5,5]$. We choose $\eta_0=0.5$ as the initial guess and plot the reconstruction at different optimization iteration steps in Figure~\ref{fig:decay3}. The reconstruction recovers the ground-truth reflection coefficient at iteration $3000$.}

\newpage
\appendix

\section{Convexity of the map $\mathcal{M}(\eta)$:}
Although we cannot demonstrate that $\mathcal{M}_\eta$ is a convex map in the entire $\eta$ function space, we can show that pointwise for every discrete $\omega$ point, the map is convex. In particular, we have the following theorem.

\begin{theorem}\label{thm:convex}
Recall the forward map to be
\begin{equation}
\mathcal{M}_\eta: {\phi} \to \Delta T(t)
\end{equation}
as defined in~\eqref{map}, with $\Delta T$ defined in~\eqref{eqn:L} where $g$ solves~\eqref{eqn:linear_PTE2} with incoming boundary condition being $\phi$. Suppose $\eta_1\geq \eta_2$ pointwise in $\omega$, then for any $\alpha,\beta\in[0,1]$ with $\alpha +\beta = 1$, one has:
\[
[\alpha\mathcal{M}_{\eta_1}+\beta\mathcal{M}_{\eta_2}] (\phi) \ge \mathcal{M}_{\alpha\eta_1+\beta\eta_2}(\phi)\,.
\]
\end{theorem}
{This suggests that if we view $\mathcal{M}$ as a function of $\eta$ as a vector and suppose it is second-order differentiable, then pointwise for $\omega$, it is convex, meaning $\partial_{ii}\mathcal{M}\geq 0$, thus the Hessian is positive along its diagonal.}
\begin{proof}
Let $g_1$, $g_2$, and $g$ solve ~\eqref{eqn:linear_PTE2} with reflection coefficients $\eta_1$, $\eta_2$ and $\alpha\eta_1+\beta\eta_2$ respectively. We essentially need to check the validity of:
\[
\Delta T_{g}(t,x=0)\leq \alpha\Delta T_{g_1}(t,x=0)+\beta\Delta T_{g_2}(t,x=0)\,,
\]
or equivalently:
\[
\langle \omega g(x=0,\cdot)\rangle\leq \alpha\langle\omega g_1(x=0,\cdot)\rangle+\beta\langle\omega g_2(x=0,\cdot)\rangle\,.
\]
In fact we will prove a stronger result that states:
\[
\alpha g_1 + \beta g_2 > g\,,
\]
for all $x,t,\mu,\omega$. To do so, we notice that with the same incoming boundary condition $\phi$, $g_i$ satisfies:
\begin{equation}
\begin{aligned}
\begin{cases}
\partial_tg_1+\mu \omega\partial_xg_1=\mathcal{L}{g_1}\\
g_1(x=0,\mu>0,\cdot) = \phi \\
g_1(x=1,\mu<0,\cdot)=\eta_1(\omega)g_1(x=1,-\mu,\cdot)
\end{cases}\,,\\
\begin{cases}
\partial_tg_2+\mu \omega\partial_xg_2=\mathcal{L}{g_2}\\
g_2(x=0,\mu>0,\cdot) = \phi \\
g_2(x=1,\mu<0,\cdot)=\eta_2(\omega)g_2(x=1,-\mu,\cdot)
\end{cases}\,
\end{aligned}
\end{equation}
We simply multiply the two PDEs with $\alpha$ and $\beta$ respectively and add them up. This gives
\begin{equation}
\begin{cases}
\partial_th+\mu \omega\partial_xh=\mathcal{L}{h}\\
h(x=0,\mu>0,\cdot) = \phi \\
h(x=1,\mu<0,\cdot)=(\alpha\eta_1+\beta\eta_2)(\omega)h(x=1,-\mu,\cdot)+\alpha\beta(\eta_1-\eta_2)(g_1-g_2)(x=1,-\mu,\cdot)
\end{cases}
\end{equation}
where $h=\alpha g_1+\beta g_2$, and we used the fact that $\alpha +\beta = 1$. On the other hand,
\begin{equation}
\begin{cases}
\partial_tg+\mu \omega\partial_xg=\mathcal{L}{g}\\
g(x=0,\mu>0,\cdot) = \phi \\
g(x=1,\mu<0,\cdot)=(\alpha\eta_1+\beta\eta_2)(\omega)g(x=1,-\mu,\cdot)
\end{cases}\,.
\end{equation}
Subtracting the PDE for $h$ and $g$, and call $H = h-g$, we have:
\begin{equation}\label{eqn:H}
\begin{cases}
\partial_tH+\mu \omega\partial_xH=\mathcal{L}{H}\\
H(x=0,\mu>0,\cdot) = 0 \\
H(x=1,\mu<0,\cdot)=(\alpha\eta_1+\beta\eta_2)(\omega)H(x=1,-\mu,\cdot)+\alpha\beta(\eta_1-\eta_2)(g_1-g_2)(x=1,-\mu,\cdot)
\end{cases}
\end{equation}
This is a phonon transport equation with zero incoming data and but a source at $x=1$. So the positivity is determined by the positivity of the source term. Since $\eta_1>\eta_2$ pointwise in $\omega$, we need to show $g_1-g_2>0$.

To show that, let $G=g_1-g_2$. Subtracting the PDEs for $g_1$ and $g_2$ gives 

\begin{equation}
\begin{cases}
\partial_tG+\mu \omega\partial_xG=\mathcal{L}{G}\\
G(x=0,\mu>0,\cdot) = 0 \\
G(x=1,\mu<0,\cdot)=\eta_2(\omega)G(x=1,-\mu,\cdot)+(\eta_1-\eta_2)(\omega)g_1(x=1,-\mu,\cdot)
\end{cases}
\end{equation}

Since both $\eta_1-\eta_2>0$ and $g_1>0$, we have $G>0$ over the entire domain, which when plugged back in~\eqref{eqn:H} indicates $H>0$, meaning,
\[
\alpha g_1+\beta g_2>g\,,
\]
concluding the theorem.
\end{proof}

\medskip

\newpage

\bibliographystyle{siamplain}
\bibliography{reference}

\end{document}